\documentclass[a4paper,12pt]{article}

\usepackage{amsmath, amsthm, amssymb}
\usepackage{bbm}

\DeclareMathOperator{\argmax}{arg\ max}
\usepackage{enumerate}
\usepackage{graphicx}
\usepackage[numbers]{natbib}
\usepackage[top=1in, bottom=1in, left=1in, right=1in]{geometry}

\title{Approximations of the Sum of States by Laplace's Method for a System of Particles with a Finite Number of Energy Levels and Application to Limit Theorems}
\author{Tomasz M. \L api\'nski \small{(84tomek@gmail.com)}\\ \\Faculty of Applied Physics and Mathematics,\\Gda\'nsk University of Technology,
Faculty of 
\\{\small ul. Gabriela Narutowicza 11/12, 80-233 Gda\'nsk, Poland}}

\begin{document}

\newtheorem{remark}{Remark}
\newtheorem{preposition}{Preposition}
\newtheorem{definition}{Definition}
\newtheorem{proposition}{Proposition}
\newtheorem{theorem}{Theorem}
\newtheorem{lemma}{Lemma}
\newtheorem{proofoflemma}{Proof of Lemma}
\newtheorem{proofofproposition}{Proof of Proposition}
\newtheorem{proofoflpreposition}{Proof of Preposition}

\maketitle

\begin{abstract}
	We consider a generic system composed of a fixed number of particles distributed over a finite number of energy levels. We make only general assumptions about system's properties and the entropy. System's constraints other than fixed number of particles can be included by appropriate reduction of system's state space. For the entropy we consider three generic cases. It can have a maximum in the interior of system's state space or on the boundary. On the boundary we can have another two cases. There the entropy can increase linearly with increase of the number of particles and in the another case grows slower than linearly. The main results are approximations of system's sum of states using Laplace's method. Estimates of the error terms are also included. As an application, we prove the law of large numbers which yields the most probable state of the system. This state is the one with the maximal entropy. We also find limiting laws for the fluctuations. These laws are different for the considered cases of the entropy. They can be mixtures of Normal, Exponential and Discrete distributions. Explicit rates of convergence are provided for all the theorems.
\end{abstract}

\begin{keywords}
	System of particles, Entropy, sum of states, Laplace's method, Law of large numbers, Central limit theorem
\end{keywords}

\begin{msc}
	41A60, 41A63, 60F05, 82B20
\end{msc}
\\
\\
The published version of this manuscript is made available at\\
https://doi.org/10.1007/s11040-020-9330-8 .
\\

\maketitle


\section{Introduction}
	\qquad We consider a system which is composed of $N$ particles distributed over $m+1$ energy levels. The energies of the levels are given by the vector $\varepsilon=(\varepsilon_{1},\varepsilon_{2},\ldots,\varepsilon_{m},\varepsilon_{m+1})$. A single system's state is represented by the vector of the occupation numbers of the levels, denoted by $(N_{1},N_{2},\ldots,N_{m},N_{m+1})$. Since the number of particles in the system is fixed, the occupation numbers satisfy
	\begin{equation}
		\label{particle_system_ensamble_constraint}
		N=\sum_{i=1}^{m+1}N_{i}.
	\end{equation}
	\indent For each accessible state, let us introduce vectors of weights of the first $m$ occupation numbers, that is, $x=(x_{1},\ldots,x_{m}):=\big(\frac{N_{1}}{N},\ldots,\frac{N_{m}}{N}\big)$. Due to constraint (\ref{particle_system_ensamble_constraint}), the last occupation number, i.e. $N_{m+1}$, is determined by the first $m$ numbers. So, the vector $x$ together with the number of particles $N$ uniquely determines a single state represented by the vector $(N_{1},\ldots,N_{m+1})$.\\
	\indent Now, let us define a set $\mathcal{D}\subset\mathbb{R}^{m}$ composed of vectors $x\in\mathbb{R}^{m}$ such that the following constraints are valid
	\begin{align*}
		&x_{1}+x_{2}+\ldots+x_{m}\leq 1,\\
		&x_{i}\geq 0,\ \text{for}\ i=1,\ldots,m,
	\end{align*}
	and also define a set $L_{N}:=\{\frac{\mathcal{N}}{N},\mathcal{N}\in\mathbb{N}^{m}\}$ for $N\geq N_{0}$, $N_{0}\in\mathbb{Z}_{+}$. We assume that $0\in\mathbb{N}^{m}$. Then, the state space of the system is represented by the set $\mathcal{D}\cap L_{N}$.\\
	\indent Initially, we assumed the system is constrained only by the number of particles. Constraints on the other properties of the system can also be included. Additional constraints can reduce the set $\mathcal{D}\cap L_{N}$ even further to a set $\mathcal{A}\cap L_{N}$, where $\mathcal{A}\subset\mathcal{D}$. So, let $\mathcal{A}\cap L_{N}$ be the set of states of the system which are accessible under given generic constraints.\\
	\indent An example of additional constraint is bounded maximal average energy per particle, i.e. $E\geq\sum_{i=1}^{m+1}\varepsilon_{i}N_{i}/N$. When such exemplary system has a large number of particles, majority of its states are in the range $(NE-\Delta,NE)$, for some small $\Delta>0$. Therefore, such system could be considered as a microcanonical ensemble, i.e., a fixed number of particles and total energy.\\
	\indent Let us consider system's entropy $S$ defined on the set $\mathcal{A}$. We assume the function $S$ is sufficiently regular, has a unique maximum and can be represented as a product of two functions, that is
	\begin{equation}
		\label{introduction_entropy}
		S(x,N)=h(N)f(x,N),
	\end{equation}
	where $f, h$ are sufficiently regular, and $h$ is also an increasing function.\\
	\indent Particular cases of the considered generic system might include features such as energy level degeneracy and indistinguishability of the particles. Such features might be reflected in the specific form of the functions $f$ and $h$, as shows the example in the end of introduction.\\
	\indent The main results of this paper are approximations of the sum of states $\Sigma(N)$ defined by
	\begin{equation}
		\label{introduction_sum_of_states}
		\Sigma(N):=\sum_{\mathcal{A}\cap L_{N}}g(x)e^{h(N)f(x,N)},
	\end{equation}
	where the function $g$ is sufficiently regular.\\
	\indent First, we develop a result for two energy levels, i.e. $m=1$,  with the maximum of the entropy on the boundary of the domain of summation. The methodology of the proof is based on the analogous result for Laplace's integral in \cite{laplace_paper}. Although, this univariate case is rather insignificant in the physics context, we need a specific estimates for the further development.\\
	\indent Our main concern is with a finite number of energy levels. We prove the results for two cases of the function $f$ in (\ref{introduction_sum_of_states}). The function $f(\cdot,N)$ can have a unique maximum in the interior of $\mathcal{A}$ or a unique non-critical maximum on the boundary $\{x:x_{1}=0\}$. For the integral instead of the sum, analogous results are proved in \cite{laplace_paper} and in a simplified form also in \cite{Laplace_method_approach_kolokoltsov}. Here we use the same methodology as in \cite{laplace_paper}, and also include an explicit remainder estimate.\\
	\indent An alternative way of approximating the sum of states for a similar class of systems is developed in \cite{Sum_of_states_alternative}.\\
	\indent Then we consider a discrete random vector $X(N):=(X_{1},X_{2},\ldots,X_{m})$, where $X(N)\in\mathcal{A}\cap L_{N}$.  We use the fundamental postulate of statistical mechanics that system's microstates are equally probable, see e.g. \cite{statistical_physics_pathria_book} and \cite{statistical_physics_reif_book}, and define the pmf of $X(N)$ to be
	\begin{equation}
		\label{particle_system_pmf}
		P(X(N)=x):=\frac{e^{S(x,N)}}{\sum_{\mathcal{A}\cap L_{N}}e^{S(y,N)}},
	\end{equation}
	where the entropy $S(x,N)$ is given by (\ref{introduction_entropy}).\\
	\indent As an application, we use the approximations of $\Sigma(N)$ to explicitely calculate the limit of $X(N)$ as $N\to\infty$. For that, we prove the law of large numbers, which can be interpreted as finding the most probable state of the system with a very large number of particles. This state is the point of maximum of the entropy (\ref{introduction_entropy}). Our next results yields the distributions of the fluctuations from the most probable state. They are different for two cases of the entropy maximum. When the maximum is in the interior of the domain, the fluctuations have Normal distribution. When the maximum is on the boundary, there can be further two cases depending on the function $h(N)$ in (\ref{introduction_entropy}). If $h(N)=N$, then the fluctuations distribution is Exponential in the direction orthogonal to the boundary of the state space and Normal in other directions. When $\lim_{N\to\infty}\frac{h(N)}{N}=0$, the fluctuations distribution is Discrete in direction orthogonal to the boundary and Normal in other directions. Explicit rate of convergence is provided for all the limit theorems.\\
	\indent  Analogous limit theorems for the integrals instead of the sums are proved in \cite{laplace_paper}. For the integral with Gaussian measure, law of large numbers and central limit theorem are proved in \cite{Laplace_method_statistical_mechanics_limit_theorems_II}, \cite{Laplace_method_statistical_mechanics_limit_theorems} and \cite{Laplace_method_statistical_mechanics}. Another application of Laplace's method to prove the limit theorems is presented in \cite{Laplace_method_weak_convergence}. \\
	\indent The results presented in this paper are applied in proofs of limit theorems in \cite{unification_paper}. There the considered system is more specific. It consist of particles that are non-interacting and indistinguishable with the average energy per particle smaller or equal to some prescribed value. Furthermore, that system has degenerate energy levels, and the number of degeneracy $G$ depends on the number of particles, that is, $G=G(N)$. Moreover, three cases of the degeneracy function $G$ are considered. The functions $f$, $h$ and the point of maximum of the entropy were derived from the system's properties and are different for each case of $G(N)$
	\begin{enumerate}[1)]
	\item $\lim_{N\to\infty}\frac{G(N)}{N}=\infty$,
		\begin{align*}
			&h(N)=N,&\\
			&f(x)=\sum_{i=1}^{m}x_{i}\ln\frac{g_{i}}{x_{i}},&\\
			&x_{i}^{*}=\frac{g_{i}}{e^{\lambda\varepsilon_{i}+\nu}},\ i=1,\ldots,m,&
		\end{align*}
	\item  $\frac{G(N)}{N}=\alpha+\beta(N),\ \text{where}\ \alpha>0, \lim_{N\to\infty}\beta(N)=0$,
		\begin{align*}
			&h(N)=N,&\\
			&f(x)=\sum_{i=1}^{m}\big((x_{i}+g_{i}\alpha)\ln(x_{i}+g_{i}\alpha)-x_{i}\ln x_{i}\big),&\\
			&x_{i}^{*}=\frac{g_{i}\alpha}{e^{\lambda\varepsilon_{i}+\nu}-1},\ i=1,\ldots,m,&
		\end{align*}
	\item $\lim_{N\to\infty}\frac{G(N)}{N}=0$,
		\begin{align*}
			&h(N)=G(N),&\\
			&f(x)=\sum_{i=1}^{m}g_{i}\ln x_{i},&\\
			&x_{i}^{*}=\frac{g_{i}}{\lambda\varepsilon_{i}+\nu},\ i=1,\ldots,m,&
		\end{align*}
	\end{enumerate}
	where $\lambda,\nu$ are some constants obtained from the systems constraints. Then with use of the theorems developed in this work author proves that for the large enough system the points of maximum given above are the most probable states. Those states are well known Maxwell-Boltzmann statistics, Bose-Einstein statistics and Zipf-Mandelbrot Law, respectively. The fluctuations are also provided.


\section{Approximation with Laplace's Method}
	\indent First, let us make several technical assumptions about the sum of states (\ref{introduction_sum_of_states}) and the functions $f,h$ and $g$. We consider a closed ball $ B_{\varepsilon}\subset\mathcal{A}$ with the center at the origin, radius $\varepsilon>0$ and volume $| B_{\varepsilon}|$. Then we specify the function $f:\mathcal{A}\times\mathbb{Z}_{+} \to \mathbb{R}$. The derivatives of $f(\cdot,N)$ up to third order exists on $ B_{\varepsilon}$ and are uniformly bounded. For all $N\geq N_{0}$, where $N_{0}\in\mathbb{Z}_{+}$, the function $f(\cdot,N)$ have a unique maximum at $x^{*}(N)\in B_{\varepsilon}$ such that
	\begin{equation}
		\label{sum_of_states_delta}
		\Delta:=\inf_{N\geq N_{0}, x\in\mathcal{A}\backslash B_{\varepsilon}}\{f(x^{*}(N),N)-f(x,N)\}>0.
	\end{equation}
	We choose the origin of our coordinate system to be the point $x^{*}=\lim_{N\to\infty}x^{*}(N)$. Further, we specify that the function $h:\mathbb{R}_{+}\to\mathbb{R}$ is positive, increasing and $\lim_{N\to\infty}\frac{h(N)}{N}=0$ or $h(N)=N$. We also specify the function $g:\mathcal{A}\to\mathbb{R}$. Assume $g$ is differentiable in $ B_{\varepsilon}$ and define constant
	\begin{align}
		\label{sum_of_states_constants_G_G1}
		&G:=\sup_{x\in\mathcal{A}}|g(x)|<\infty,\quad G^{(1)}:=\sup_{x\in B_{\varepsilon}}\|Dg(x)\|<\infty,\\
		\label{sum_of_states_bound}
		&C>0,\ C\geq\sum_{\mathcal{A}\cap L_{N}}\frac{1}{N^{m}}e^{h(N_{0})f(x,N)},\ \text{for all}\ N\geq N_{0},
	\end{align}
	where $\|.\|$ is a max norm and $D$ is a differential operator, that is, $D:=(\frac{\partial}{\partial x_{1}},\frac{\partial}{\partial x_{2}},\ldots,\frac{\partial}{\partial x_{m}})$.
	Let us assume that the sum (\ref{introduction_sum_of_states}) is finite and specify the two cases of its function $f$
	\begin{enumerate}[(a)]
		\item $f(\cdot,N)$ has a nondegenerate maximum in the interior of $ B_{\varepsilon}$ and introduce constant
			\begin{align}
				\label{sum_of_states_interior_constant_F'2}
				&F'^{(2)}:=\inf_{x\in B_{\varepsilon}, N\geq N_{0}}\|D^{2}f(x,N)^{-\frac{1}{2}}\|^{-2}>0,\\
				\label{sum_of_states_interior_constant_F'2det}
				&F'^{(2)}_{det}:=\inf_{x\in B_{\varepsilon}, N\geq N_{0}}\sqrt{|\det D^{2}f(x,N)|}>0,\\
				\label{sum_of_states_interior_constant_F2}
				&F^{(2)}:=\sup_{x\in B_{\varepsilon},N\geq N_{0}}\|D^{2}f(x,N)\|<\infty,\\
				\label{sum_of_states_interior_constant_F3}
				&F^{(3)}:=\sup_{x\in B_{\varepsilon},N\geq N_{0}}\|D^{3}f(x,N)\|<\infty,
			\end{align}
		\item $f(\cdot,N)$ has a unique maximum on the boundary $\{x:x_{1}=0\}$ and $\frac{\partial f(x^{*}(N),N)}{\partial x_{1}}<0$. We also introduce constant
			\begin{align}
				\label{sum_of_states_boundary_constant_F'1}
				&F'^{(1)}:=\inf_{x\in B_{\varepsilon}, N\geq N_{0}}\bigg|\frac{\partial f(x,N)}{\partial x_{1}}\bigg|>0,\\
				\label{sum_of_states_boundary_constant_F'2}
				&F'^{(2)}:=\inf_{x\in B_{\varepsilon}, N\geq N_{0}}\|D_{y}^{2}f(x,N)^{-\frac{1}{2}}\|^{-2}>0,\\
				\label{sum_of_states_boundary_constant_F'2det}
				&F'^{(2)}_{det}:=\inf_{x\in B_{\varepsilon}, N\geq N_{0}}\sqrt{|\det D_{y}^{2}f(x,N)|}>0,\\
				\label{sum_of_states_boundary_constant_F2}
				&F^{(2)}:=\sup_{x\in B_{\varepsilon},N\geq N_{0}}\|D^{2}f(x,N)\|<\infty,\\
				\label{sum_of_states_boundary_constant_F3}
				&F^{(3)}:=\sup_{x\in B_{\varepsilon},N\geq N_{0}}\|D^{3}f(x,N)\|<\infty,
			\end{align}
			where $y=(x_{2},\ldots,x_{m})$ and $D_{y}$ is a differential operator in that coordinates.\\
			Furthermore, we assume that on every section $ B_{\varepsilon}(x_{1})=\{ y:(x_{1},y)\in B_{\varepsilon}\}, x_{1}\in[0,\varepsilon)$ we have a unique nondegenerate maximum of f.
	\end{enumerate}

	\begin{remark}
		The situation when the boundary of the domain is $\{x:x_{1}=a\}$ with $a\in\mathbb{Q}_{+}$ can be reduced to the case of the boundary $\{x:x_{1}=0\}$, if we only consider $N$ such that $Na\in\mathbb{Z}$. This is because for those values, the structure of lattice $L_{N}$ is preserved after appropriate shift of the coordinate system.
	\end{remark}


\subsection{One Dimensional Entropy}
	For the case (b) of the function $f$ and $\mathcal{A}=[0,\infty)$ we define a set
	\begin{equation*}
		U_{N}:=\bigg\{x:x\leq\frac{1}{h(N)^{1-\delta}},x\in\mathcal{A}\bigg\},
	\end{equation*}
	where the parameter $\delta\in(0,1)$. Additionally, let us define $N_{1}:=\max\big\{\big\lceil h^{-1}\big(\varepsilon^{\frac{1}{\delta-1}}\big)\big\rceil,N_{0}\big\}$. For all $N\geq N_{1}$, the set $U_{N}\subset B_{\varepsilon}$ and we have the following results

	\begin{proposition}
		\label{univariate_preposition}
		For the case (b) of the function $f$ and $\mathcal{A}=[0,\infty)$, the following approximation holds
		\begin{equation*}
			\sum_{U_{N}\cap L_{N}}e^{h(N)f'(0,N)x}=\frac{1}{1-\exp\big(\frac{h(N)}{N}f'(0,N)\big)}\Big(1+\omega_{UB}(N)\exp\big(-|f'(0,N)|h(N)^{\delta}\big)\Big),
		\end{equation*}
		where $\omega_{UB}(N)=O(1)\ \text{as}\ N\to\infty$ and
		\begin{align*}
			|\omega_{UB}&(N)|<\frac{N}{h(N)}\frac{1-\exp\big(\frac{h(N)}{N}f'(0,N)\big)}{|f'(0,N)|}\exp\bigg(\frac{h(N)}{N}|f'(0,N)|\bigg).
		\end{align*}
	\end{proposition}

	\begin{theorem}
		\label{univariate_theorem}
		For the case (b) of the function $f$ and $\mathcal{A}=[0,\infty)$, the following approximation holds
		\begin{align*}
			\sum_{\mathcal{A}\cap L_{N}}g(x)e^{h(N)f(x,N)}=e^{h(N)f(0,N)}\frac{1}{1-\exp\big(\frac{h(N)}{N}f'(0,N)\big)}\Bigg(g(0)+\frac{\omega_{B}(N)}{h(N)^{1-2\delta}}\Bigg),
		\end{align*}
		where $\delta\in(0,\frac{1}{2})$, $\omega_{B}(N)=O(1)\ \text{as}\ N\to\infty$ and 
		\begin{align*}
			|&\omega_{B}(N)|\leq\bigg(\frac{GF^{(2)}}{2}+G^{(1)}h(N)^{-\delta}\bigg)\exp\bigg(\frac{F^{(2)}}{2h(N)^{1-\delta}}\bigg)\Big(1+\omega_{UB}(N)\exp\big(-F'^{(1)}h(N)^{\delta}\big)\Big)+\\
			&\omega_{UB}(N)Gh(N)^{1-2\delta}\exp\big(-F'^{(1)}h(N)^{\delta}\big)+Gh(N)^{1-2\delta}\bigg(1-\exp\bigg(\frac{h(N)}{N}f'(0,N)\bigg)\bigg)\times\\
			&\Big(\exp\big(-F'^{(1)}h(N)^{\delta}\big)(N\varepsilon+1)+NC\exp\big(-h(N)\Delta-h(N_{0})(f(0,N)-\Delta)\big)\Big).
		\end{align*}
		where $\omega_{UB}$ is defined in Proposition \ref{univariate_preposition}.
	\end{theorem}

	\begin{proof}[Proof of Proposition \ref{univariate_preposition}]
		Let $x=\frac{i}{N}$, $i=0,\ldots,I_{N}$, with $I_{N}=\big\lfloor\frac{N}{h(N)^{1-\delta}}\big\rfloor$. Then we have 
		\begin{equation*}
				\sum_{U_{N}\cap L_{N}}e^{h(N)f'(0,N)x}=\sum_{i=0}^{I_{N}}\exp\bigg(\frac{h(N)}{N}f'(0,N)i\bigg),
		\end{equation*}
		which is equal
		\begin{align*}
			\sum_{U_{N}\cap L_{N}}e^{h(N)f'(0,N)x}&=\sum_{i=0}^{\infty}\exp\bigg(\frac{h(N)}{N}f'(0,N)i\bigg)-\sum_{i>I_{N}}\exp\bigg(\frac{h(N)}{N}f'(0,N)i\bigg)=\\
			&\frac{1}{1-\exp\big(\frac{h(N)}{N}f'(0,N)\big)}-\sum_{i>I_{N}}\exp\bigg(\frac{h(N)}{N}f'(0,N)i\bigg),
		\end{align*}
		by the formula for the summation of the geometric series. Then, we estimate the last term by the simple approximation of the sum with an integral
		\begin{align*}
			&\sum_{i>I_{N}}\exp\bigg(\frac{h(N)}{N}f'(0,N)i\bigg)\leq\int_{I_{N}}^{\infty}\exp\bigg(\frac{h(N)}{N}f'(0,N)y\bigg)dy<\\
			&\frac{N}{h(N)|f'(0,N)|}\exp\bigg(-|f'(0,N)|\bigg(h(N)^{\delta}-\frac{h(N)}{N}\bigg)\bigg)=\\
			&\frac{N}{h(N)|f'(0,N)|}\exp\bigg(\frac{h(N)}{N}|f'(0,N)|\bigg)\exp\big(-|f'(0,N)|h(N)^{\delta}\big),
		\end{align*}
		where $I_{N}>\frac{N}{h(N)^{1-\delta}}-1$ and the expression in the exponent is negative due to (\ref{sum_of_states_boundary_constant_F'1}). Hence we get the result of the Proposition.
	\end{proof}

	\begin{proof}[Proof of Theorem \ref{univariate_theorem}]
		Let us introduce $S_{B}(N)$, $\Sigma_{B}(N)$ and using the Taylor's Theorem decompose $\Sigma(N)$  
		\begin{align*}
			S_{B}(N)&:=g(0)e^{h(N)f(0,N)}\frac{1}{1-\exp\big(\frac{h(N)}{N}f'(0,N)\big)},\\
			\Sigma_{B}(N)&:=g(0)e^{h(N)f(0,N)}\sum_{U_{N}\cap L_{N}}\exp\big(h(N)f'(0,N)x\big),\\
			\Sigma(N)&=\Sigma_{11}(N)+\Sigma_{12}(N)+\Sigma_{2}(N)+\Sigma_{3}(N):=g(0)\sum_{U_{N}\cap L_{N}}e^{h(N)f(x,N)}+\\
			&\sum_{U_{N}\cap L_{N}}g'(x_{\theta}(N))xe^{h(N)f(x,N)}+\sum_{( B_{\varepsilon}\backslash U_{N})\cap L_{N}}g(x)e^{h(N)f(x,N)}+\sum_{(\mathcal{A}\backslash B_{\varepsilon})  \cap L_{N}}g(x)e^{h(N)f(x,N)}.
		\end{align*}
		Here and everywhere in the proofs $x_{\theta}$ denotes a point between $x$ and the point of the expansion. It might be different in a different instances. Now, we put together the above expressions
		\begin{equation*}
			|\Sigma(N)-S_{B}(N)|\leq|\Sigma_{11}(N)-\Sigma_{B}(N)|+|\Sigma_{B}(N)-S_{B}(N)|+|\Sigma_{12}(N)|+|\Sigma_{2}(N)|+|\Sigma_{3}(N)|.
		\end{equation*}
		For $|\Sigma_{11}(N)-\Sigma_{B}(N)|$ we use the second order Taylor's Theorem to obtain			
		\begin{align}
			\label{Extended_Laplace_proof_laplace_approximation_one_dimension_main_term}
			|\Sigma_{11}&(N)-\Sigma_{B}(N)|\leq|g(0)|\sum_{U_{N}\cap L_{N}}\exp\big(h(N)f(0,N)+h(N)f'(0,N)x\big)\times\notag\\
			&\bigg|\exp\bigg(\frac{1}{2}h(N)f''(x_{\theta}(N),N)x^{2}\bigg)-1\bigg|.
		\end{align}
		The second term in the Taylor's Theorem can be bounded, that is $|f''(x_{\theta},N)x^{2}|\leq F^{(2)}x^{2}$, where $F^{(2)}$ is defined by (\ref{sum_of_states_boundary_constant_F2}). Next, using result (\ref{Extended_Laplace_proof_laplace_approximation_one_dimension_main_term}) with the inequality $|e^{t}-1|\leq|t|e^{|t|}$ and the fact that for any $x\in U_{N}$, $x\leq h(N)^{-1+\delta}$ yields
		\begin{align}
			|\Sigma_{11}&(N)-\Sigma_{B}(N)|\leq\frac{1}{2}GF^{(2)}e^{h(N)f(0,N)}h(N)^{-1+2\delta}\sum_{U_{N}\cap L_{N}}\exp\big(h(N)f'(0,N)x\big)\times\notag\\
			\label{Extended_Laplace_proof_laplace_approximation_one_dimension_main_term_2}
			&\exp\bigg(\frac{1}{2}F^{(2)}h(N)^{-1+2\delta}\bigg),
		\end{align}
		where $G$ is defined by (\ref{sum_of_states_constants_G_G1}). We need the last term in the inequality (\ref{Extended_Laplace_proof_laplace_approximation_one_dimension_main_term_2}) to be bounded as $N\to\infty$, hence we set $\delta\in\big(0,\frac{1}{2}\big)$. Then, with use of Proposition \ref{univariate_preposition} we obtain the estimate
		\begin{align*}
			\label{Extended_Laplace_proof_laplace_approximation_m_dimension_interior_of_domain_approximation_of_integral_I11_IG1}
			|\Sigma_{11}(N)-\Sigma_{B}(N)|&\leq\frac{1}{2}GF^{(2)}\exp\bigg(\frac{F^{(2)}}{2h(N)^{1-2\delta}}\bigg)e^{h(N)f(0,N)}h(N)^{-1+2\delta}\times\\
			&\frac{1}{1-\exp\big(\frac{h(N)}{N}f'(0,N)\big)}\Big(1+\omega_{UB}(N)\exp\big(-|f'(a,N)|h(N)^{\delta}\Big).
		\end{align*} 
		Next expression to approximate, $|\Sigma_{B}(N)-S_{B}(N)|$, can be directly obtained from Proposition \ref{univariate_preposition}
		\begin{equation*}
			|\Sigma_{B}(N)-S_{B}(N)|\leq \omega_{UB}(N)|g(0)|e^{h(N)f(0,N)}\frac{1}{1-\exp\big(\frac{h(N)}{N}f'(0,N)\big)}\exp\big(-|f'(0,N)|h(N)^{\delta}\big).
		\end{equation*} 
		Now, let us consider the sum $\Sigma_{12}(N)$. Here again, we apply the second order Taylor's Theorem and since $g$ has a bounded derivative in $U_{N}$ we obtain 
		\begin{equation*}
			|\Sigma_{12}(N)|\leq G^{(1)}e^{h(N)f(0,N)}h(N)^{-1+\delta}\sum_{U_{N}\cap L_{N}}e^{h(N)f'(0,N)x}\exp\bigg(\frac{1}{2}h(N)F^{(2)}x^{2}\bigg),
		\end{equation*}
		where $G^{(1)}$ is defined by (\ref{sum_of_states_constants_G_G1}). Further, applying Proposition \ref{univariate_preposition} and using that $x\leq h(N)^{-1+\delta}$ we obtain
		\begin{align*}
			|\Sigma_{12}(N)|&\leq G^{(1)}\exp\bigg(\frac{F^{(2)}}{2h(N)^{1-2\delta}}\bigg)e^{h(N)f(0,N)}h(N)^{-1+\delta}\frac{1}{1-\exp\big(\frac{h(N)}{N}f'(0,N)\big)}\times\\
			&\Big(1+\omega_{UB}(N)\exp\big(-|f'(0,N)|h(N)^{\delta}\big)\Big).
		\end{align*}	
		For $|\Sigma_{2}(N)|$ we apply the first order Taylor's Theorem which yields
		\begin{equation*}
			\label{Laplace_integral_interior_proof_f_estimate}
			f(x,N)\leq f(x^{*}(N),N)-F'^{(1)}x,
		\end{equation*} 
		where $F'^{(1)}$ is defined by (\ref{sum_of_states_boundary_constant_F'1}). Then we substitute it into $|\Sigma_{2}(N)|$ and get
		\begin{equation*}
			|\Sigma_{2}(N)|\leq e^{h(N)f(0,N)}\sum_{( B_{\varepsilon}\backslash U_{N})\cap L_{N}}|g(x)|\exp\big(-h(N)F'^{(1)}x\big).
		\end{equation*}
		Since in the set $ B_{\varepsilon}\backslash U_{N}$ function $g$ is bounded by $G$ and $x>\frac{1}{h(N)^{1-\delta}}$, hence
		\begin{align*}
 			|\Sigma_{2}(N)|\leq Ge^{h(N)f(0,N)}\exp\big(-F'^{(1)}h(N)^{\delta}\big)\sum_{ B_{\varepsilon}\cap L_{N}}1.
		\end{align*}
		The number of elements in the set is bounded by $\varepsilon N+1$. Therefore
		\begin{align*}
 			|\Sigma_{2}(N)|\leq Ge^{h(N)f(0,N)}\exp\big(-h(N)^{\delta}F'^{(1)}\big)(\varepsilon N+1).
		\end{align*}
		In case of $|\Sigma_{3}(N)|$ we have the following upper bound
		\begin{align*}
			|\Sigma_{3}(N)|&\leq e^{h(N)f(0,N)}\sum_{(\mathcal{A}\backslash B_{\varepsilon})\cap L_{N}}|g(x)|\exp\big(h(N_{0})(f(x,N)-f(0,N))-(h(N)-h(N_{0}))\Delta\big)\leq\\
			&\exp\big((h(N)-h(N_{0}))(f(0,N)-\Delta)\big)\sum_{(\mathcal{A}\backslash B_{\varepsilon})\cap L_{N}}|g(x)|e^{h(N_{0})f(x,N)}\leq\\
			&GCe^{h(N)f(0,N)}N\exp\big(-h(N)\Delta-h(N_{0})(f(0,N)-\Delta)\big),
		\end{align*}
		where the last inequality is due to assumption (\ref{sum_of_states_bound}).\\
		Then, we combine the above approximations 
		\begin{align*}
			\label{Laplaca_approximation_m-dimensional_interior_proof_combined_estimate}
			|&\Sigma(N)-S_{B}(N)|\leq e^{h(N)f(0,N)}h(N)^{-1+2\delta}\frac{1}{1-\exp\big(\frac{h(N)}{N}f'(0,N)\big)}\Bigg[\bigg(\frac{GF^{(2)}}{2}+G^{(1)}h(N)^{-\delta}\bigg)\times\\
			&\exp\bigg(\frac{F^{(2)}}{2h(N)^{1-\delta}}\bigg)\Big(1+\omega_{UB}(N)\exp\big(-F'^{(1)}N^{\delta}\big)\Big)+\omega_{UB}(N)Gh(N)^{1-2\delta}\exp\big(-F'^{(1)}h(N)^{\delta}\big)+\\
			&Gh(N)^{1-2\delta}\bigg(\exp\big(-F'^{(1)}h(N)^{\delta}\big)(\varepsilon N+1)+CN\exp\big(-h(N)\Delta-h(N_{0})(f(0,N)-\Delta)\big)\bigg)\times\\
			&\bigg(1-\exp\bigg(\frac{h(N)}{N}f'(0,N)\bigg)\bigg)\Bigg].
		\end{align*}
	\end{proof}


\subsection{Multivariate Entropy with Maximum in the Interior}
	For the case (a) of the function $f$ we define sets
	\begin{align*}
		&U_{N}:=\bigg\{x:|x-x^{*}(N)|\leq\frac{1}{h(N)^{\frac{1}{2}-\delta}}\bigg\},\\
		&V_{N,y}:=\bigg\{x:y_{i}-\frac{1}{2N}\leq x_{i}<y_{i}+\frac{1}{2N}, i\in1,2,\ldots,m\bigg\},\\
		&V_{N}:=\{x:x\in V_{N,y},y\in U_{N}\cap L_{N}\},
	\end{align*}
	where the parameter $\delta\in\big(0,\frac{1}{2}\big)$. Further, let us define $N_{1}:=\max\big\{\big\lceil h^{-1}\big(\varepsilon^{\frac{2}{2\delta-1}}\big)\big\rceil,N_{0}\big\}$. For all $N\geq N_{1}$, the set $U_{N}\subset B_{\varepsilon}$ we have the following results 
	\begin{proposition}
		\label{multivariate_interior_preposition}
		For the case (a) of the function $f$, the following approximation holds
		\begin{align*}
			&\sum_{U_{N}\cap L_{N}}\exp\bigg(\frac{1}{2}h(N)(x-x^{*}(N))^{T}D^{2}f(x^{*}(N),N)(x-x^{*}(N))\bigg)=\\
			&N^{m}\bigg(\frac{2\pi}{h(N)}\bigg)^{\frac{m}{2}}\bigg(\frac{1}{\sqrt{|\det D^{2}f(x^{*}(N),N)|}}+\omega_{UI}(N)\frac{h(N)^{1/2+(m+1)\delta}}{N}\bigg),
		\end{align*}
		where $\delta\in\big(0,\frac{1}{2(m+1)}\big)$, $\omega_{UI}(N)=O(1)\ \text{as}\ N\to\infty$ and
		\begin{align*}
			|&\omega_{UI}(N)|\leq\frac{2^{\frac{m}{2}}F^{(2)}}{\Gamma(\frac{m}{2}+1)}\bigg(1+\frac{\sqrt{m}h(N)^{1/2-\delta}}{N}\bigg)^{m+1}+\frac{N}{h(N)^{1/2+\delta}}\exp\bigg(-\frac{1}{2}h(N)^{2\delta}F'^{(2)}\bigg)\times\\
			&\Bigg[\frac{2^{-\frac{m}{2}+1}}{\Gamma\big(\frac{m}{2}\big)}\exp\bigg(\frac{\sqrt{m}F'^{(2)}h(N)^{1/2+\delta}}{N}+\frac{mF'^{(2)}h(N)}{2N^{2}}\bigg)\Bigg(\bigg(1+\frac{\sqrt{m}h(N)^{1/2-\delta}}{N}\bigg)^{m}-\\
			&\bigg(1-\frac{\sqrt{m}h(N)^{1/2-\delta}}{N}\bigg)^{m}\Bigg)+\frac{1}{F'^{(2)}_{det}}\exp\bigg(\frac{1}{2}h(N_{0})^{2\delta}F'^{(2)}\bigg)\frac{1}{h(N)^{m\delta}}\bigg(\frac{h(N)}{h(N_{0})}\bigg)^{\frac{m}{2}}\Bigg].
		\end{align*}
	\end{proposition}

	\begin{theorem}
		\label{multivariate_interior_theorem}
		For the case (a) of the function $f$, the following approximation holds
		\begin{align*}
			\sum_{\mathcal{A}\cap L_{N}}&g(x)e^{h(N)f(x,N)}=e^{h(N)f(x^{*}(N),N)}N^{m}\bigg(\frac{2\pi}{h(N)}\bigg)^{\frac{m}{2}}\times\\
			&\Bigg(\frac{g(x^{*}(N))}{\sqrt{|\det D^{2}f(x^{*}(N),N)|}}+\omega_{I}(N)\frac{1}{h(N)^{1/2-3\delta}}+\omega_{UI}(N)G\frac{h(N)^{1/2+(m+1)\delta}}{N}\Bigg),
		\end{align*}
		where $\delta\in\big(0,\frac{1}{3(m+1)}\big)$, $\omega_{I}(N)=O(1)\ \text{as}\ N\to\infty$ and 
		\begin{align*}
			|\omega_{I}&(N)|\leq\bigg(\frac{1}{\sqrt{|\det D^{2}f(x^{*}(N),N)|}}+\omega_{UI}(N)\frac{h(N)^{1/2+(m+1)\delta}}{N}\bigg)\exp\bigg(\frac{F^{(3)}}{6h(N)^{1/2-3\delta}}\bigg)\times\\
			&\bigg(\frac{GF^{(3)}}{6}+G^{(1)}h(N)^{-2\delta}\bigg)+G\bigg(\frac{2\pi}{h(N)}\bigg)^{-\frac{m}{2}}h(N)^{1/2-3\delta}\Bigg(\frac{\pi^{\frac{m}{2}}}{\Gamma(\frac{m}{2}+1)}\exp\bigg(-\frac{1}{2}h(N)^{2\delta}F'^{(2)}\bigg)\times\\
			&\bigg(\varepsilon+\frac{\sqrt{m}}{N}\bigg)^{m}+C\exp\big(-h(N)\Delta-h(N_{0})(f(x^{*},N)-\Delta)\big)\Bigg),
		\end{align*}
		where $\omega_{UI}$ is defined in Proposition \ref{multivariate_interior_preposition}.
	\end{theorem}

	\begin{proof}[Proof of Proposition \ref{multivariate_interior_preposition}]
		Let us define $I(N)$ and $\Sigma'(N)$ 
		\begin{align*}
			I(N)&:=N^{m}\bigg(\frac{2\pi}{h(N)}\bigg)^{\frac{m}{2}}\frac{1}{\sqrt{|\det D^{2}f(x^{*}(N),N)|}}=\\
			&N^{m}\int_{\mathbb{R}^{m}}\exp\bigg(\frac{1}{2}h(N)(x-x^{*}(N))^{T}D^{2}f(x^{*}(N),N)(x-x^{*}(N))\bigg)dx,\\
			\Sigma'(N)&:=\sum_{U_{N}\cap L_{N}}\exp\bigg(\frac{1}{2}h(N)(x-x^{*}(N))^{T}D^{2}f(x^{*}(N),N)(x-x^{*}(N))\bigg),
		\end{align*} 
		and decompose $I$ into four integrals with the indicator function
		\begin{align*}
			I(N)&=I_{1}(N)+I_{2}(N)-I_{3}(N)+I_{4}(N):=\\
			&N^{m}\int_{\mathbb{R}^{m}}\bigg(\mathbbm{1}_{V_{N}}(x)+\mathbbm{1}_{U_{N}\backslash V_{N}}(x)-\mathbbm{1}_{V_{N}\backslash U_{N}}(x)+\mathbbm{1}_{\mathbb{R}^{m}\backslash U_{N}}(x)\bigg)\times\\
			&\exp\bigg(\frac{1}{2}h(N)(x-x^{*}(N))^{T}D^{2}f(x^{*}(N),N)(x-x^{*}(N))\bigg)dx.
		\end{align*} 
		Then decompose $I_{1}(N)$ into a smaller integrals and use the Taylor's Theorem 
		\begin{align*}
			&I_{1}(N)=I_{11}(N)+I_{12}(N):=\\
			&N^{m}\sum_{y\in U_{N}\cap L_{N}}\int_{V_{N,y}}\exp\bigg(\frac{1}{2}h(N)(y-x^{*}(N))^{T}D^{2}f(x^{*}(N),N)(y-x^{*}(N))\bigg)dx+N^{m}\times\\
			&\sum_{y\in U_{N}\cap L_{N}}\int_{V_{N,y}}(x-y)^{T}D\exp\bigg(\frac{1}{2}h(N)(x-x^{*}(N))^{T}D^{2}f(x^{*}(N),N)(x-x^{*}(N))\bigg)\Big\vert_{x=x_{\theta}(y)}dx.
		\end{align*}  
		We combine above decompositions into	
		\begin{equation*}
			|\Sigma'(N)-I(N)|\leq|\Sigma'(N)-I_{11}(N)|+|I_{12}(N)|+|I_{2}(N)-I_{3}(N)|+|I_{4}(N)|,
		\end{equation*}
		and approximate each term separately.\\
		\indent First, we consider the expression $|\Sigma'(N)-I_{11}(N)|$. For the integral $I_{11}(N)$, we have			
		\begin{align*}
			I_{11}(N)&=N^{m}\sum_{y\in U_{N}\cap L_{N}}\exp\bigg(\frac{1}{2}h(N)(y-x^{*}(N))^{T}D^{2}f(x^{*}(N),N)(y-x^{*}(N))\bigg)\int_{V_{N,y}}1dx=\\
			&\sum_{y\in U_{N}\cap L_{N}}\exp\bigg(\frac{1}{2}h(N)(y-x^{*}(N))^{T}D^{2}f(x^{*}(N),N)(y-x^{*}(N))\bigg),
		\end{align*}
		as the integral is equal to the volume of the hypercube $V_{N,y}$, that is, $N^{-m}$. Hence $|\Sigma'(N)-I_{11}(N)|=0.$\\
		\indent Then we approximate $|I_{12}(N)|$. Since for any $x\in V_{N,y}$, $|x-y|\leq \frac{\sqrt{m}}{N}$, hence
		\begin{align}
			\label{Extended_Laplace_proof_laplace_approximation_finite_dimension_interior_I12}
			|I&_{12}(N)|\leq N^{m}\sum_{y\in U_{N}\cap L_{N}}\frac{\sqrt{m}}{N}\times\\
			&\int_{V_{N,y}}\bigg\vert D\exp\bigg(\frac{1}{2}h(N)(x-x^{*}(N))^{T}D^{2}f(x^{*}(N),N)(x-x^{*}(N))\bigg)\Big\vert_{x=x_{\theta}(y)}\bigg\vert dx.\notag
		\end{align}
		Next, we estimate the derivative in (\ref{Extended_Laplace_proof_laplace_approximation_finite_dimension_interior_I12})
		\begin{align*}
			&\bigg\vert D\exp\bigg(\frac{1}{2}h(N)(x-x^{*}(N))^{T}D^{2}f(x^{*}(N),N)(x-x^{*}(N))\bigg)\Big\vert_{x=x_{\theta}(y)}\bigg\vert\leq \\
			&\frac{1}{2}h(N)\max_{i}\bigg|2\sum_{j=1}^{m}\frac{\partial^{2} f(x^{*}(N),N)}{\partial x_{i}\partial x_{j}}(x_{j}-x_{j}^{*}(N))\bigg|\times\\
			&\exp\bigg(\frac{1}{2}h(N)(x_{\theta}(y)-x^{*}(N))^{T}D^{2}f(x^{*}(N),N)(x_{\theta}(y)-x^{*}(N))\bigg)\leq \\
			&h(N)mF^{(2)}|x-x^{*}(N)| \exp\bigg(\frac{1}{2}h(N)(x_{\theta}(y)-x^{*}(N))^{T}D^{2}f(x^{*}(N),N)(x_{\theta}(y)-x^{*}(N))\bigg),
		\end{align*}
		where $F^{(2)}$ is defined by (\ref{sum_of_states_interior_constant_F2}). Since in the integration we included the points outside $U_{N}$ and within $V_{N}$, hence we have 
		\begin{equation*}
			|x-x^{*}(N)|\leq \frac{1}{h(N)^{1/2-\delta}}+\frac{\sqrt{m}}{N}=\frac{1}{h(N)^{1/2-\delta}}\bigg(1+\frac{\sqrt{m}h(N)^{1/2-\delta}}{N}\bigg).
		\end{equation*}
		Using that and the fact that the volume of $V_{N,y}$ is $N^{-m}$ we obtain
		\begin{align}
			\label{Extended_Laplace_proof_m_dimensional_I12}
			|I_{12}(N)&|\leq m^{3/2}F^{(2)}\frac{h(N)^{1/2+\delta}}{N}\bigg(1+\frac{\sqrt{m}h(N)^{1/2-\delta}}{N}\bigg)\times\notag\\
			&\sum_{y\in U_{N}\cap L_{N}}\exp\bigg(\frac{1}{2}h(N)(x_{\theta}(y)-x^{*}(N))^{T}D^{2}f(x^{*}(N),N)(x_{\theta}(y)-x^{*}(N))\bigg)\leq\notag\\
			&m^{3/2}F^{(2)}\frac{h(N)^{1/2+\delta}}{N}\bigg(1+\frac{\sqrt{m}h(N)^{1/2-\delta}}{N}\bigg)\sum_{U_{N}\cap L_{N}}1,
		\end{align}
		where the last inequality is because $D^{2}f$ is negative definite in $U_{N}\subset B_{\varepsilon}$, hence occurring exponent can be bounded by $1$. Now we estimate the size of the sum $\sum_{U_{N}\cap L_{N}}1$. It is clear that this sum is bounded by the number of the hypercubes $V_{N,y}, y\in L_{N}$ that intersects $U_{N}$. The sphere of the radius $h(N)^{-1/2+\delta}+\sqrt{m}N^{-1}$ and dimension $m$ contains all such hypercubes. Therefore, this sphere's volume divided the volume of hypercube $V_{N,y}$, that is
		\begin{equation*}
			\bigg(\frac{N}{h(N)^{1/2-\delta}}+\sqrt{m}\bigg)^{m}\frac{\pi^{m/2}}{\Gamma\big(\frac{m}{2}+1\big)},
		\end{equation*}
		is an upper bound for the number of $V_{N,y}$ that intersects $U_{N}$. Putting that into the estimate (\ref{Extended_Laplace_proof_m_dimensional_I12}) yields
		\begin{equation*}
			|I_{12}(N)|\leq\frac{m^{3/2}F^{(2)}}{\Gamma(\frac{m}{2}+1)}N^{m}\bigg(\frac{\pi}{h(N)}\bigg)^{\frac{m}{2}}\frac{h(N)^{1/2+(m+1)\delta}}{N}\bigg(1+\frac{\sqrt{m}h(N)^{1/2-\delta}}{N}\bigg)^{m+1}.
		\end{equation*}
		For the approximation of $|I_{2}(N)-I_{3}(N)|$ let us introduce a set 
		\begin{equation*}
		\widetilde{U}_{N}:=\bigg\{x:\frac{1}{h(N)^{1/2-\delta}}-\frac{\sqrt{m}}{N}\leq |x-x^{*}(N)|\leq\frac{1}{h(N)^{1/2-\delta}}+\frac{\sqrt{m}}{N}\bigg\}.
		\end{equation*}
		Since $\widetilde{U}_{N}$ contains the domains of the integration of $I_{2}(N)$ and $I_{3}(N)$ we have 
		\begin{align}
			\label{Extended_Laplace_proof_laplace_approximation_finite_dimension_I2-I3}
			|I_{2}&(N)-I_{3}(N)|\leq N^{m}\int_{\widetilde{U}_{N}}\exp\bigg(\frac{1}{2}h(N)(x-x^{*}(N))^{T}D^{2}f(x^{*}(N),N)(x-x^{*}(N))\bigg)dx\leq\notag\\
			&N^{m}\int_{\widetilde{U}_{N}}\exp\bigg(-\frac{1}{2}h(N)F'^{(2)}|x-x^{*}(N)|^{2}\bigg)dx\leq\\
			&N^{m}\exp\bigg(-\frac{1}{2}h(N)^{2\delta}F'^{(2)}\bigg|1-\frac{\sqrt{m}h(N)^{1/2-\delta}}{N}\bigg|^{2}\bigg)\int_{\widetilde{U}_{N}}1dx,\notag
		\end{align}
		where we used the fact that
		\begin{align*}
			&(x-x^{*}(N))^{T}D^{2}f(x^{*}(N),N)(x-x^{*}(N))\leq -\|D^{2}f(x^{*}(N),N)\| |x-x^{*}(N)|^{2}\leq \\
			&-F'^{(2)}|x-x^{*}(N)|^{2},
		\end{align*}
		with $F'^{(2)}$ defined by (\ref{sum_of_states_interior_constant_F'2}). Next, we calculate the integral in (\ref{Extended_Laplace_proof_laplace_approximation_finite_dimension_I2-I3}) by performing the change of coordinates system to the spherical
		\begin{align*}
			&\int_{\widetilde{U}_{N}}1dx=\frac{2\pi^{m/2}}{m\Gamma(\frac{m}{2})}\bigg[\bigg(\frac{1}{h(N)^{1/2-\delta}}+\frac{\sqrt{m}}{N}\bigg)^{m}-\bigg(\frac{1}{h(N)^{1/2-\delta}}-\frac{\sqrt{m}}{N}\bigg)^{m}\bigg].
		\end{align*}
		Hence, we obtain the following estimate for $|I_{G2}(N)-I_{G3}(N)|$
		\begin{align*}
			&|I_{2}(N)-I_{3}(N)|\leq\frac{2}{\Gamma(\frac{m}{2})}N^{m}\bigg(\frac{\pi}{h(N)}\bigg)^{\frac{m}{2}}h(N)^{m\delta}\exp\bigg(-\frac{1}{2}h(N)^{2\delta}F'^{(2)}\bigg|1-\frac{\sqrt{m}h(N)^{1/2-\delta}}{N}\bigg|^{2}\bigg)\times\\
			&\bigg[\bigg(1+\frac{\sqrt{m}h(N)^{1/2-\delta}}{N}\bigg)^{m}-\bigg(1-\frac{\sqrt{m}h(N)^{1/2-\delta}}{N}\bigg)^{m}\bigg].
		\end{align*}
		Finally, we approximate the last integral $I_{4}(N)$
		\begin{align*}
 			|I_{4}(N)|\leq& N^{m}\int_{\mathbb{R}^{m}\backslash U_{N}}\exp\bigg(\frac{1}{2}(h(N)-h(N_{0}))(x-x^{*}(N))^{T}D^{2}f(x^{*}(N),N)(x-x^{*}(N))\bigg)\times\\
			&\exp\bigg(\frac{1}{2}h(N_{0})(x-x^{*}(N))^{T}D^{2}f(x^{*}(N),N)(x-x^{*}(N))\bigg)dx\leq\\
			&\exp\big(-\frac{1}{2}h(N)^{2\delta}F'^{(2)}+\frac{1}{2}h(N_{0})^{2\delta}F^{(2)}\bigg)\times\\
			&\int_{\mathbb{R}^{m}}\exp\bigg(\frac{1}{2}h(N_{0})(x-x^{*}(N))^{T}D^{2}f(x^{*}(N),N)(x-x^{*}(N))\bigg)dx\leq \\
			&\frac{1}{F'^{(2)}_{det}}\bigg(\frac{2\pi}{h(N_{0})}\bigg)^{\frac{m}{2}}N^{m}\exp\bigg(-\frac{1}{2}h(N)^{2\delta}F'^{(2)}+\frac{1}{2}h(N_{0})^{2\delta}F'^{(2)}\bigg),
		\end{align*}
		where $F'^{(2)}_{det}$ is given by (\ref{sum_of_states_interior_constant_F'2det}).\\	
		Then we combine all approximations 
		\begin{align*}
			\label{Laplaca_approximation_m-dimensional_interior_proof_combined_estimate}
			|\Sigma'(&N)-I(N)|\leq N^{m}\bigg(\frac{2\pi}{h(N)}\bigg)^{\frac{m}{2}}\frac{h(N)^{1/2+(m+1)\delta}}{N}\Bigg[\frac{2^{-\frac{m}{2}}F^{(2)}m^{\frac{3}{2}}}{\Gamma(\frac{m}{2}+1)}\bigg(1+\frac{\sqrt{m}h(N)^{1/2-\delta}}{N}\bigg)^{m+1}+\\
			&\frac{2^{-\frac{m}{2}+1}}{\Gamma(\frac{m}{2})}\frac{N}{h(N)^{1/2+\delta}}\exp\bigg(-\frac{1}{2}h(N)^{2\delta}F'^{(2)}\bigg|1-\frac{\sqrt{m}h(N)^{1/2-\delta}}{N}\bigg|^{2}\bigg)\times\\
			&\Bigg(\bigg(1+\frac{\sqrt{m}h(N)^{1/2-\delta}}{N}\bigg)^{m}-\bigg(1-\frac{\sqrt{m}h(N)^{1/2-\delta}}{N}\bigg)^{m}\Bigg)+\\
			&\frac{1}{F'^{(2)}_{det}}\frac{N}{h(N)^{1/2+(m+1)\delta}}\bigg(\frac{h(N)}{h(N_{0})}\bigg)^{\frac{m}{2}}\exp\bigg(-\frac{1}{2}h(N)^{2\delta}F'^{(2)}+\frac{1}{2}h(N_{0})^{2\delta}F'^{(2)}\bigg)\Bigg],
		\end{align*}
		and for the error term to decrease as $N\to\infty$, we set $\delta\in\big(0,\frac{1}{2(m+1)}\big)$.
	\end{proof}
	
	\begin{proof}[Proof of Theorem \ref{multivariate_interior_theorem}]
		Let us introduce $I_{G}(N)$, $\Sigma_{G}(N)$ and using the Taylor's Theorem decompose $\Sigma(N)$  
		\begin{align*}
			I_{G}&(N):=N^{m}g(x^{*}(N))e^{h(N)f(x^{*}(N),N)}\bigg(\frac{2\pi}{h(N)}\bigg)^{\frac{m}{2}}\frac{1}{\sqrt{|\det D^{2}f(x^{*}(N),N)|}}=\\
			&N^{m}g(x^{*}(N))e^{h(N)f(x^{*}(N),N)}\times\\
			&\int_{\mathbb{R}^{m}}\exp\bigg(\frac{1}{2}h(N)(x-x^{*}(N))^{T}D^{2}f(x^{*}(N),N)(x-x^{*}(N))\bigg)dx,\\
			\Sigma_{G}&(N):=g(x^{*}(N))e^{h(N)f(x^{*}(N),N)}\times\\
			&\sum_{U_{N}\cap L_{N}}\exp\bigg(\frac{1}{2}h(N)(x-x^{*}(N))^{T}D^{2}f(x^{*}(N),N)(x-x^{*}(N))\bigg),\\
			\Sigma&(N)=\Sigma_{11}(N)+\Sigma_{12}(N)+\Sigma_{2}(N)+\Sigma_{3}(N):=g(x^{*}(N))\sum_{U_{N}\cap L_{N}}e^{h(N)f(x,N)}+\\
			&\sum_{U_{N}\cap L_{N}}Dg(x_{\theta}(N))^{T}(x-x^{*}(N))e^{h(N)f(x,N)}+\sum_{( B_{\varepsilon}\backslash U_{N})\cap L_{N}}g(x)e^{h(N)f(x,N)}+\\
			&\sum_{(\mathcal{A}\backslash B_{\varepsilon})\cap L_{N}}g(x)e^{h(N)f(x,N)}.
		\end{align*} 
		We combine the above decompositions into	
		\begin{equation*}
			|\Sigma(N)-I_{G}(N)|\leq|\Sigma_{11}(N)-\Sigma_{G}(N)|+|\Sigma_{G}(N)-I_{G}(N)|+|\Sigma_{12}(N)|+|\Sigma_{2}(N)|+|\Sigma_{3}(N)|.
		\end{equation*}
		For $|\Sigma_{11}(N)-\Sigma_{G}(N)|$ we use the third order Taylor's Theorem to obtain			
		\begin{align}
			\label{Laplace_Approximation_m_dimensional_interior_theorem_proof_Sigma11}
			&|\Sigma_{11}(N)-\Sigma_{G}(N)|=|g(x^{*}(N))|\times\notag\\
			&\sum_{U_{N}}\exp\bigg(h(N)f(x^{*}(N),N)+\frac{1}{2}h(N)(x-x^{*}(N))^{T}D^{2}f(x^{*}(N),N)(x-x^{*}(N))\bigg)\times\\
			&\bigg[\exp\bigg(\frac{1}{6}h(N)D^{3}f(x_{\theta}(N),N)(x-x^{*}(N))^{3}\bigg)-1\bigg],\notag
		\end{align}
		due to $Df(x^{*}(N),N)^{T}(x-x^{*}(N))=0$, since $x^{*}(N)$ is critical point. The third term in the Taylor's Theorem can be bounded
		\begin{equation*}
			|D^{3}f(x_{\theta})x^{3}|\leq\|D^{3}f(x_{\theta})\| |x|^{3}\leq F^{(3)}|x|^{3},
		\end{equation*}
		where $F^{(3)}$ is defined by (\ref{sum_of_states_interior_constant_F3}). Next, using the result (\ref{Laplace_Approximation_m_dimensional_interior_theorem_proof_Sigma11}) with the inequality $|e^{t}-1|\leq|t|e^{|t|}$ and the fact that for $x\in U_{N}$, $|x-x^{*}(N)|\leq h(N)^{-1/2+\delta}$ yields
		\begin{align*}
			|\Sigma_{11}&(N)-\Sigma_{G}(N)|\leq\frac{1}{6}GF^{(3)}e^{h(N)f(x^{*}(N),N)}h(N)^{-1/2+3\delta}\times\\
			&\sum_{U_{N}\cap L_{N}}\exp\bigg(\frac{1}{2}h(N)(x-x^{*}(N))^{T}D^{2}f(x^{*}(N),N)(x-x^{*}(N))\bigg)\exp\bigg(\frac{1}{6}F^{(3)}h(N)^{-1/2+3\delta}\bigg),
		\end{align*}
		with $G$ defined by (\ref{sum_of_states_constants_G_G1}). In order to bound the last term in the above estimate, we set $\delta\in\big(0,\frac{1}{6}\big)$. Then with use of Proposition \ref{multivariate_interior_preposition} we obtain the following estimate
		\begin{align*}
			\label{Extended_Laplace_proof_laplace_approximation_m_dimension_interior_of_domain_approximation_of_integral_I11_IG1}
			|\Sigma_{11}(N)-\Sigma_{G}(N)|&\leq\frac{1}{6}GF^{(3)}\exp\bigg(\frac{F^{(3)}}{6h(N)^{1/2-3\delta}}\bigg)e^{h(N)f(x^{*}(N),N)}N^{m}\bigg(\frac{2\pi}{h(N)}\bigg)^{\frac{m}{2}}h(N)^{-1/2+3\delta}\times\\
			&\bigg(\frac{1}{\sqrt{|\det D^{2}f(x^{*}(N),N)|}}+\omega_{UI}(N)\frac{h(N)^{1/2+(m+1)\delta}}{N}\bigg),
		\end{align*} 
		with $\delta\in\big(0,\frac{1}{2(m+1)}\big)$. For the estimate to be valid for all $m\in\mathbb{Z}_{+}$, we set $\delta\in\big(0,\frac{1}{3(m+1)}\big)$.\\
		\indent Next expression to approximate, $|\Sigma_{G}(N)-I_{G}(N)|$, can be directly obtained from Proposition \ref{multivariate_interior_preposition}
		\begin{equation*}
			|\Sigma_{G}(N)-I_{G}(N)|\leq \omega_{UI}(N)g(x^{*}(N))e^{h(N)f(x^{*}(N),N)}N^{m}\bigg(\frac{2\pi}{h(N)}\bigg)^{\frac{m}{2}}\omega_{UI}(N)\frac{h(N)^{1/2+(m+1)\delta}}{N}.
		\end{equation*} 
		Now, let us consider the sum $\Sigma_{12}(N)$. Here again, we apply the third order Taylor's Theorem and since in $U_{N}$ the derivative of $g$ is bounded by $G^{(1)}$, we obtain 
		\begin{align*}
			|\Sigma_{12}&(N)|\leq G^{(1)}e^{h(N)f(x^{*}(N),N)}h(N)^{-1/2+\delta}\exp\bigg(\frac{1}{6}h(N)F^{(3)}|x-x^{*}(N)|^{3}\bigg)\times\\
			&\sum_{U_{N}\cap L_{N}}\exp\bigg(\frac{1}{2}h(N)(x-x^{*}(N))^{T}D^{2}f(x^{*}(N),N)(x-x^{*}(N))\bigg),
		\end{align*}
		where the constant $G^{(1)}$ is defined by (\ref{sum_of_states_constants_G_G1}). Further, applying Proposition \ref{multivariate_interior_preposition} and using that $|x-x^{*}(N)|\leq h(N)^{-1/2+\delta}$ in $U_{N}$ we obtain
		\begin{align*}
			|\Sigma_{12}&(N)|\leq G^{(1)}\exp\bigg(\frac{F^{(3)}}{6h(N)^{1/2-3\delta}}\bigg)e^{h(N)f(x^{*}(N),N)}N^{m}\bigg(\frac{2\pi}{h(N)}\bigg)^{\frac{m}{2}}h(N)^{-1/2+\delta}\times\\
			&\bigg(\frac{1}{\sqrt{|\det D^{2}f(x^{*}(N),N)|}}+\omega_{UI}(N)\frac{h(N)^{1/2+(m+1)\delta}}{N}\bigg).
		\end{align*}	
		For $|\Sigma_{2}(N)|$, we apply the second order Taylor's Theorem to $f$ to get
		\begin{equation*}
			\label{Laplace_integral_interior_proof_f_estimate}
			f(x,N)\leq f(x^{*}(N),N)-\frac{1}{2}F'^{(2)}|x-x^{*}(N)|^{2},
		\end{equation*} 
		and then estimate $|\Sigma_{2}(N)|$
		\begin{equation*}
			|\Sigma_{2}(N)|\leq e^{h(N)f(x^{*}(N),N)}\sum_{( B_{\varepsilon}\backslash U_{N})\cap L_{N}}g(x)\exp\bigg(-\frac{1}{2}NF'^{(2)}|x-x^{*}(N)|^{2}\bigg).
		\end{equation*}
		Since in the set $ B_{\varepsilon}\backslash U_{N}$, the function $g$ is bounded by $G$ and $|x-x^{*}(N)|>\frac{1}{N^{1/2-\delta}}$, hence
		\begin{equation}
			\label{Laplace_approximation_m_dimensional_interior_theorem_proof_Sigma2}
 			|\Sigma_{2}(N)|\leq Ge^{h(N)f(x^{*}(N),N)}\exp\bigg(-\frac{1}{2}h(N)^{2\delta}F'^{(2)}\bigg)\sum_{ B_{\varepsilon}\cap L_{N}}1.
		\end{equation}
		Now, we estimate the size of the sum $\sum_{U_{N}\cap L_{N}}1$. It is clear that this sum is bounded by the number of the hypercubes $V_{N,y}$, $y\in L_{N}$ that intersects $B_{\varepsilon}$. The sphere of the radius $\varepsilon+\sqrt{m}N^{-1}$ and dimension $m$ contains all such hypercubes. Therefore, this sphere volume divided by the volume of $V_{N,y}$, that is 
		\begin{equation*}
			\frac{\pi^{m/2}}{\Gamma(\frac{m}{2}+1)}(N\varepsilon+\sqrt{m})^{m},
		\end{equation*}
		is an upperbound for the number of $V_{N,y}$ that intersects $B_{\varepsilon}$. Adding that to the estimate (\ref{Laplace_approximation_m_dimensional_interior_theorem_proof_Sigma2}) yields
		\begin{align*}
 			|\Sigma_{2}(N)|\leq \frac{G\pi^{m/2}}{\Gamma(\frac{m}{2}+1)} N^{m}e^{h(N)f(x^{*}(N),N)}\exp\bigg(-\frac{1}{2}h(N)^{2\delta}F'^{(2)}\bigg)\bigg(\varepsilon+\frac{\sqrt{m}}{N}\bigg)^{m}.
		\end{align*}
		In case of $|\Sigma_{3}(N)|$ we have the following upper bound
		\begin{align*}
			|\Sigma_{3}(N)|&\leq e^{h(N)f(x^{*}(N),N)}\times\\
			&\sum_{(\mathcal{A}\backslash B_{\varepsilon})\cap L_{N}}|g(x)|\exp\big(h(N_{0})(f(x,N)-f(x^{*}(N),N))-(h(N)-h(N_{0}))\Delta\big)\leq\\
			&\exp\big((h(N)-h(N_{0}))(f(x^{*}(N),N)-\Delta)\big)\sum_{(\mathcal{A}\backslash B_{\varepsilon})\cap L_{N}}|g(x)|e^{h(N_{0})f(x,N)}\leq\\
			&GCe^{h(N)f(x^{*}(N),N)}N^{m}\exp\big(-h(N)\Delta-h(N_{0})(f(x^{*}(N),N)-\Delta)\big),
		\end{align*}
		where the last inequality is due to assumptions (\ref{sum_of_states_constants_G_G1}) and (\ref{sum_of_states_bound}).\\
		Then we combine the above approximations 
		\begin{align*}
			\label{Laplaca_approximation_m-dimensional_interior_proof_combined_estimate}
			|&\Sigma(N)-I_{G}(N)|\leq h(N)^{-1/2+3\delta}e^{h(N)f(x^{*}(N),N)}N^{m}\bigg(\frac{2\pi}{h(N)}\bigg)^{\frac{m}{2}}\times\\
			&\Bigg[\bigg(\frac{1}{\sqrt{|\det D^{2}f(x^{*}(N),N)|}}+\omega_{UI}(N)\frac{h(N)^{1/2+(m+1)\delta}}{N}\bigg)\exp\bigg(\frac{F^{(3)}}{6h(N)^{1/2-3\delta}}\bigg)\times\\
			&\bigg(\frac{GF^{(3)}}{6}+G^{(1)}h(N)^{-2\delta}\bigg)+G\bigg(\frac{2\pi}{h(N)}\bigg)^{-\frac{m}{2}}h(N)^{1/2-3\delta}\bigg(\frac{\pi^{\frac{m}{2}}}{\Gamma(\frac{m}{2}+1)}\exp\bigg(-\frac{1}{2}h(N)^{2\delta}F'^{(2)}\bigg)\times\\
			&\bigg(\varepsilon+\frac{\sqrt{m}}{N}\bigg)^{m}+C\exp\big(-h(N)\Delta-h(N_{0})(f(x^{*},N)-\Delta)\big)\bigg)\Bigg]+\\
			&\omega_{UI}(N)g(x^{*}(N))\omega_{UI}(N)e^{h(N)f(x^{*}(N),N)}N^{m}\bigg(\frac{2\pi}{h(N)}\bigg)^{\frac{m}{2}}\frac{h(N)^{1/2+(m+1)\delta}}{N}.
		\end{align*}
	\end{proof}


\subsection{Multivariate Entropy with Maximum on the Boundary}
	
	\begin{theorem}
		\label{multivariate_boundary_theorem}
		For the case (b) of the function $f$, the following approximation holds
		\begin{align*}
			&\sum_{\mathcal{A}\cap L_{N}\cap\{x:x_{1}\geq 0\}}g(x)e^{h(N)f(x,N)}=e^{N f(x^{*}(N),N)}N^{m-1}\bigg(\frac{2\pi}{h(N)}\bigg)^{\frac{m-1}{2}}\frac{1}{1-\exp\big(\frac{h(N)}{N}\frac{\partial f(x^{*}(N),N)}{\partial x_{1}}\big)}\times\\
			&\Bigg(\frac{g(x^{*}(N))}{\sqrt{\big|\det D_{y}^{2}f(x^{*}(N),N)\big|}}+\omega_{1}(N)\frac{1}{h(N)^{1/2-3\delta}}+\omega_{2}(N)\frac{h(N)^{1/2+(m+1)\delta}}{N}\Bigg),
		\end{align*}
		valid for all $N\geq N_{1}$ with $N_{1}:=\max\big\{\big\lceil h^{-1}\big(\varepsilon^{\frac{1}{\delta-1}}\big)\big\rceil,\big\lceil h^{-1}\big(\varepsilon^{\frac{2}{2\delta-1}}\big)\big\rceil,N_{0}\big\}$, $\delta=\big(0,\frac{1}{3(m+1)}\big)$, $\omega_{1}(N)=O(1)$, $\omega_{2}(N)=O(1)$, $\omega_{C}(N)=O(1)$ as $N\to\infty$ and
		\begin{align*}
			&\omega_{1}(N)=\omega_{I}(N)+\omega_{I}(N)\frac{\omega_{B2}(N)}{h(N)^{1-2\delta}}+\omega_{B1}(N)\frac{1}{F'^{(2)}_{det}}\frac{1}{h(N)^{1/2+\delta}}+\omega_{C}(N),\\
			&\omega_{2}(N)=\omega_{UI}(N)G+\omega_{UI}(N)\frac{G\omega_{B2}(N)}{h(N)^{1-2\delta}},
		\end{align*}
		where $\omega_{UI}$ is defined in Proposition \ref{multivariate_interior_preposition}, $\omega_{I}$ in Theorem \ref{multivariate_interior_theorem} and $\omega_{B1}$, $\omega_{B2}$, $\omega_{C}$ are 
		\begin{align*}
			|\omega_{B1}&(N)|\leq\Big(1+\omega_{UB}(N)\exp\big(-F'^{(1)}h(N)^{\delta}\big)\Big)\exp\bigg(\frac{F^{(2)}}{2h(N)^{1-\delta}}\bigg)\times\\
			&\bigg(\frac{GF^{(2)}}{2}+G^{(1)}h(N)^{-\delta}+\frac{m^{2}GF^{(3)}}{2F'^{(2)}}h(N)^{-\delta}\bigg)+Gh(N)^{1-2\delta}\exp\big(-F'^{(1)}h(N)^{\delta}\big)\times\\
			&\Bigg(\omega_{UB}(N)+(N\varepsilon+1)\bigg(1-\exp\bigg(\frac{h(N)}{N}\frac{\partial f(x^{*}(N),N)}{\partial x_{1}}\bigg)\bigg)\Bigg),\\
			|\omega_{B2}&(N)|\leq\Big(1+\omega_{UB}(N)\exp\big(-F'^{(1)}h(N)^{\delta}\big)\Big)\exp\bigg(\frac{F^{(2)}}{2h(N)^{1-\delta}}\bigg)\frac{F^{(2)}}{2}+h(N)^{1-2\delta}\times\\
			&\exp\big(-F'^{(1)}h(N)^{\delta}\big)\Bigg(\omega_{UB}(N)+(N\varepsilon+1)\bigg(1-\exp\bigg(\frac{h(N)}{N}\frac{\partial f(x^{*}(N),N)}{\partial x_{1}}\bigg)\bigg)\Bigg),\\
			|\omega_{C}&(N)|\leq N\bigg(\frac{2\pi}{h(N)}\bigg)^{-\frac{m-1}{2}}h(N)^{1/2-3\delta}\bigg(1-\exp\bigg(\frac{h(N)}{N}\frac{\partial f(x^{*}(N),N)}{\partial x_{1}}\bigg)\bigg)\times\\
			&GC\exp\big(-h(N)\Delta-h(N_{0})(f(x^{*}(N),N)-\Delta)\big),
		\end{align*}
		where $\omega_{UB}$ inside $\omega_{B1}$ and $\omega_{B2}$ is defined in Proposition \ref{univariate_preposition}.
	\end{theorem}
	
	\begin{remark}
		The  situation when the boundary is an arbitrary hyperplane with rational coefficients can be reduced to the case with boundary $\{ x:x_{1}= 0\}$. This is because after appropriate rotation of the coordinate system, the structure of the lattice, which is essential for the application of Theorem \ref{multivariate_boundary_theorem} is preserved. That is, all the points of the domain are on the equally spaced hyperplanes parallel to the boundary $\{x:x_{1}=0\}$.
	\end{remark}

	\begin{proof} 
		We decompose $\Sigma(N)$ into
		\begin{align*}
			\Sigma(&N)=\Sigma_{1}(N)+\Sigma_{2}(N):=\sum_{ B_{\varepsilon}\cap L_{N}\cap\{x:x_{1}\geq 0\}}g(x)e^{h(N)f(x,N)}dx+\\
			&\sum_{(\mathcal{A}\cap\{x:x_{1}\geq 0\}\backslash B_{\varepsilon})\cap L_{N}}g(x)e^{h(N)f(x,N)}dx,
		\end{align*}
		and approximate $\Sigma_{2}(N)$ as the sum $\Sigma_{3}(N)$ in the previous proof.\\
		Since $\mathcal{A}\backslash B_{\varepsilon}\subset\mathcal{A}\backslash U_{N}$
		\begin{align*}
			|\Sigma_{2}(N)|&\leq\exp\big((h(N)-h(N_{0}))(f(x^{*}(N),N)-\Delta)\big)\sum_{(\mathcal{A}\cap\{x:x_{1}\geq 0\}\backslash B_{\varepsilon})\cap L_{N}}|g(x)|e^{h(N_{0})f(x,N)}dx\leq\\
			&N^{m}e^{h(N)f(x^{*}(N),N)}GC\exp\big(-h(N)\Delta-h(N_{0})(f(x^{*}(N),N)-\Delta)\big).
		\end{align*}
		Then we express $\Sigma_{1}(N)$ as
		\begin{equation*}
			\Sigma_{1}(N)=\sum_{x_{1}\in B_{\varepsilon}\cap L_{N}, x_{1}\geq 0}\Sigma_{1}(x_{1},N),
		\end{equation*}
		where the sum is over the values of coordinate $x_{1}$ of the points in $ B_{\varepsilon}\cap L_{N}$. Further, we have
		\begin{equation*}
			\Sigma_{1}(x_{1},N):=\sum_{ B_{\varepsilon}\cap L_{N}(x_{1}), x_{1}\geq 0}g(x_{1},y)e^{h(N)f(x_{1},y,N)},
		\end{equation*}
		where $ B_{\varepsilon}\cap L_{N}(x_{1})=\{y:(x_{1},y)\in B_{\varepsilon}\cap L_{N}\}$. Next, we apply Theorem \ref{multivariate_interior_theorem} to $\Sigma_{1}(x_{1},N)$
		\begin{align*}
			\Sigma_{1}(x_{1},N)=&e^{h(N)f(x_{1},y^{*}(x_{1},N),N)}N^{m-1}\bigg(\frac{2\pi}{h(N)}\bigg)^{\frac{m-1}{2}}\Bigg(\frac{g(x_{1},y^{*}(x_{1},N))}{\sqrt{|\det D_{y}^{2}f(x_{1},y^{*}(x_{1},N),N)|}}+\\
			&\omega_{I}(x_{1},N)\frac{1}{h(N)^{1/2-3\delta}}+\omega_{UI}(x_{1},N)G\frac{h(N)^{1/2+(m+1)\delta}}{N}\Bigg),
		\end{align*}
		where $y^{*}(x_{1},N)=\argmax_{y\in B_{\varepsilon}(x_{1})}f(x_{1},y,N)$. As the summation is over the set $ B_{\varepsilon}\cap L_{N}(x_{1})$, the constants which occurs as a result of application of Theorem 2 can be replaced by the appropriate constants for the larger set $ B_{\varepsilon}$, which are independent of $x_{1}$, that is (\ref{sum_of_states_boundary_constant_F'2}), (\ref{sum_of_states_boundary_constant_F'2det}), (\ref{sum_of_states_boundary_constant_F2}) and (\ref{sum_of_states_boundary_constant_F3}). Then, we apply Theorem \ref{univariate_theorem} to $\Sigma_{1}(N)$
		\begin{align*}
			\Sigma_{1}(N)=&e^{N f(x^{*}(N),N)}N^{m-1}\bigg(\frac{2\pi}{h(N)}\bigg)^{\frac{m-1}{2}}\frac{1}{1-\exp\big(\frac{h(N)}{N}\frac{\partial f(x^{*}(N),N)}{\partial x_{1}}\big)}\Bigg(\frac{g(x^{*}(N))}{\sqrt{|\det D_{y}^{2}f(x^{*}(N),N)|}}+\\
			&+\frac{\omega_{B1}(N)}{F'^{(2)}_{det}h(N)^{1-2\delta}}+\frac{\omega_{I}(N)}{h(N)^{1/2-3\delta}}+\frac{\omega_{I}(N)\omega_{B2}(N)}{h(N)^{\frac{3}{2}-5\delta}}+\\
			&+\frac{\omega_{UI}(N)Gh(N)^{1/2+(m+1)\delta}}{N}+\frac{\omega_{UI}(N)\omega_{B2}(N)G}{h(N)^{1-2\delta}}\frac{h(N)^{1/2+(m+1)\delta}}{N}\Bigg),
		\end{align*}
		where  $(0,y^{*}(0,N))=x^{*}(N)$. Since Theorem \ref{univariate_theorem} was applied on the curve $y^{*}(x_{1},N)$, the constants in the estimate of $\omega_{B1}(N)$ and $\omega_{B2}(N)$ can also be replaced by the constants for the larger set $ B_{\varepsilon}$ i.e. (\ref{sum_of_states_boundary_constant_F'1}), (\ref{sum_of_states_boundary_constant_F'2}) and (\ref{sum_of_states_boundary_constant_F2}). 
		\begin{align*}
			|\omega_{B1}&(N)|\leq\Big(1+\omega_{UB}(N)\exp\big(-F'^{(1)}h(N)^{\delta}\big)\Big)\exp\bigg(\frac{F^{(2)}}{2h(N)^{1-\delta}}\bigg)\times\\
			&\bigg(\frac{GF^{(2)}}{2}+G^{(1)}h(N)^{-\delta}+\frac{m^{2}GF^{(3)}}{2F'^{(2)}}h(N)^{-\delta}\bigg)+h(N)^{1-2\delta}G\exp\big(-F'^{(1)}h(N)^{\delta}\big)\times\\
			&\Bigg(\omega_{UB}(N)+(N\varepsilon+1)\bigg(1-\exp\bigg(\frac{h(N)}{N}\frac{\partial f(x^{*}(N),N)}{\partial x_{1}}\bigg)\bigg)\Bigg).\\
			|\omega_{B2}&(N)|\leq\Big(1+\omega_{UB}(N)\exp\big(-F'^{(1)}h(N)^{\delta}\big)\Big)\exp\bigg(\frac{F^{(2)}}{2h(N)^{1-\delta}}\bigg)\frac{F^{(2)}}{2}+h(N)^{1-2\delta}\times\\
			&\exp\big(-F'^{(1)}h(N)^{\delta}\big)\Bigg(\omega_{UB}(N)+(N\varepsilon+1)\bigg(1-\exp\bigg(\frac{h(N)}{N}\frac{\partial f(x^{*}(N),N)}{\partial x_{1}}\bigg)\bigg)\Bigg).
		\end{align*}
		Then we combine the above result with the estimate of $\Sigma_{2}(N)$ to obtain the final result.
	\end{proof}


\section{Limit Theorems}

	For the function $f$ in (\ref{introduction_sum_of_states}), let us additionally assume
	\begin{equation}
		\label{probability_equation_for_f}
		f(x,N)=f(x)+\epsilon(N)\sigma(x),
	\end{equation}
	where $\sigma(x)$, $f(x)$ are functions that possess the derivatives at $x^{*}$ up to second order and the function $\epsilon(N)>0$ for all $N$, and $\epsilon(N)\to 0$ as $N\to\infty$. Furthermore, for the case (a) of $f(x,N)$ presented in the beginning of Section 2, we assume that the function $f(x)$ has a nondegenerate maximum at $x^{*}$. In the case (b), we assume $\frac{\partial f(x^{*})}{\partial x_{1}}<0$ and w.r.t. coordinates $(x_{2},\ldots,x_{m})$ function $f(x)$ has a nondegenerate maximum at the point $x^{*}$ on the boundary $\{x:x_{1}=0\}$. For both cases, equation (\ref{probability_equation_for_f}) implies
	\begin{equation}
		\label{probability_estimate_maximums}
		x^{*}(N)=x^{*}+\epsilon(N)O(1),\ N\to\infty.
	\end{equation} 
	Let us consider two cases of the pmf of $X(N)$ given by (\ref{particle_system_pmf})
	\begin{align}
		\label{pmf_interior}
		(a)&\ P(X(N)=x)=\frac{e^{h(N)f(x,N)}}{\sum_{\mathcal{A}\cap L_{N}}e^{h(N)f(y,N)}},\text{ for the case (a) of }f,\\
		\label{pmf_boundary}
		(b)&\ P(X(N)=x)=\frac{e^{h(N)f(x,N)}}{\sum_{\mathcal{A}\cap\{x:x_{1}\geq 0\}\cap L_{N}}e^{h(N)f(y,N)}},\text{ for the case (b) of }f.
	\end{align}


\subsection{Weak Law of Large Numbers}
	
	\begin{theorem}[Weak law of large numbers]
		As $N\to \infty$, the random vector $X(N)$ converges in distribution to the constant $x^{*}$ and the following estimate of the mgf holds
		\begin{equation*}
			M_{X(N)}(\xi)=e^{\xi^{T}x^{*}}\bigg(1+\frac{O(1)}{h(N)^{1/2-3\delta}}+O(1)\frac{h(N)^{1/2+(m+1)\delta}}{N}+O(1)\epsilon(N)\bigg),\ N\to\infty,
		\end{equation*}
		where $\delta\in\big(0,\frac{1}{3(m+1)}\big)$.
	\end{theorem}
	
	\begin{remark}
		For this and the following limit theorems the convergence error term can be explicitly estimated with use of the results from the Section 2.
	\end{remark}
	
	\begin{proof}
		To prove the convergence of $X(N)$, it is sufficient to prove the convergence of its moment generating function
		\begin{align*}
			(a)&M_{X(N)}(\xi)=\frac{\sum_{\mathcal{A}\cap L_{N}}e^{\xi^{T}x}e^{h(N)f(x,N)}}{\sum_{\mathcal{A}\cap L_{N}}e^{h(N)f(x,N)}},\\
			(b)&M_{X(N)}(\xi)=\frac{\sum_{\mathcal{A}\cap\{x:x_{1}\geq 0\}\cap L_{N}}e^{\xi^{T}x}e^{h(N)f(x,N)}}{\sum_{\mathcal{A}\cap\{x:x_{1}\geq 0\}\cap L_{N}}e^{h(N)f(x,N)}},
		\end{align*}
		where $|\xi|<h$, for some $h>0$. We approximate the denominator for the case (a) of $f$ with use of Theorem \ref{multivariate_interior_theorem}, and for the case (b) using Theorem \ref{multivariate_boundary_theorem}
		\begin{align*}
			(a)&\sum_{\mathcal{A}\cap L_{N}}e^{h(N)f(x,N)}=e^{h(N)f(x^{*}(N),N)}N^{m}\bigg(\frac{2\pi}{h(N)}\bigg)^{\frac{m}{2}}\frac{1}{\sqrt{|\det D^{2}f(x^{*}(N),N)|}}\times\\
			&\bigg(1+\frac{O(1)}{h(N)^{1/2-3\delta}}+O(1)\frac{h(N)^{1/2+(m+1)\delta}}{N}\bigg),\\
			(b)&\sum_{\mathcal{A}\cap\{x:x_{1}\geq 0\}\cap L_{N}}e^{h(N)f(x,N)}=e^{h(N)f(x^{*}(N),N)}N^{m-1}\bigg(\frac{2\pi}{h(N)}\bigg)^{\frac{m-1}{2}}\times\\
			&\frac{1}{\Big(1-\exp\big(\frac{h(N)}{N}\frac{\partial f(x^{*}(N),N)}{\partial x_{1}}\big)\Big)\sqrt{|\det D_{y}^{2}f(x^{*}(N),N)|}}\times\\
			&\bigg(1+\frac{O(1)}{h(N)^{1/2-3\delta}}+O(1)\frac{h(N)^{1/2+(m+1)\delta}}{N}\bigg).
		\end{align*}	
		Then we do the same for the numerator
		\begin{align*}
			(a)&\sum_{\mathcal{A}\cap L_{N}}\exp\big(\xi^{T}x+h(N)f(x,N)\big)=\exp\big(\xi^{T}x^{*}(N)+h(N)f(x^{*}(N),N)\big)N^{m}\bigg(\frac{2\pi}{h(N)}\bigg)^{\frac{m}{2}}\times\\
			&\frac{1}{\sqrt{|\det D^{2}f(x^{*}(N),N)|}}\bigg(1+\frac{O(1)}{h(N)^{1/2-3\delta}}+O(1)\frac{h(N)^{1/2+(m+1)\delta}}{N}\bigg),\\
			(b)&\sum_{\mathcal{A}\cap\{x:x_{1}\geq 0\}\cap L_{N}}\exp\big(\xi^{T}x+h(N)f(x,N)\big)=\exp\big(\xi^{T}x^{*}(N)+h(N)f(x^{*}(N),N)\big)\times\\
			&N^{m-1}\bigg(\frac{2\pi}{h(N)}\bigg)^{\frac{m-1}{2}}\frac{1}{\Big(1-\exp\big(\frac{h(N)}{N}\frac{\partial f(x^{*}(N),N)}{\partial x_{1}}\big)\Big)\sqrt{|\det D_{y}^{2}f(x^{*}(N),N)|}}\times\\
			&\bigg(1+\frac{O(1)}{h(N)^{1/2-3\delta}}+O(1)\frac{h(N)^{1/2+(m+1)\delta}}{N}\bigg).
		\end{align*}
		Dividing the approximation of denominator by the approximation of numerator yields for the both cases
		\begin{equation*}
			\label{limit_proof_limit_mgf}
			 M_{X(N)}(\xi)=\exp\big(\xi^{T}x^{*}(N)\big)\bigg(1+\frac{O(1)}{h(N)^{1/2-3\delta}}+O(1)\frac{h(N)^{1/2+(m+1)\delta}}{N}\bigg).
		\end{equation*}
 		Next, we use the estimate (\ref{probability_estimate_maximums}) and the Taylor's expansion for the exponent function to obtain
		\begin{equation*}
			 M_{X(N)}(\xi)=\exp\big(\xi^{T}x^{*}\big)\big(1+\epsilon(N)\xi^{T}O(1)\big)\bigg(1+\frac{O(1)}{h(N)^{1/2-3\delta}}+O(1)\frac{h(N)^{1/2+(m+1)\delta}}{N}\bigg),
		\end{equation*}
		which leads to the result of the theorem.
	\end{proof}


\subsection{Central Limit Theorem}	

	Here, let us consider the function $\epsilon(N)$ introduced in (\ref{probability_equation_for_f}) to be $\epsilon(N)=o\Big(\frac{1}{\sqrt{h(N)}}\Big)$ as $N\to\infty$. Then for the case (a) of $f$, we have the following result
	\begin{proposition}
		\label{probability_preposition}
		For the function $\widetilde{f}(x,N):=f(x,N)+\frac{1}{\sqrt{h(N)}}\xi^{T}(x^{*}-x)$ with $\xi>0$, the  following approximations holds
		\begin{align}
			\label{probability_tilda_estimates_maximums}
			&\widetilde{x}^{*}(N)-x^{*}(N)=D^{2}f(x^{*})^{-1}\frac{\xi}{\sqrt{h(N)}}+\frac{O(1)}{h(N)},\\
			\label{probability_tilda_estimates_functions}
			&\widetilde{f}(\widetilde{x}^{*}(N),N)-f(x^{*}(N),N)=\frac{1}{2h(N)}\xi^{T}D^{2}f(x^{*})^{-1}\xi+\frac{O(1)}{h(N)^{3/2}}+\frac{O(1)\epsilon(N)}{\sqrt{h(N)}},\\
			\label{probability_tilda_estimates_determinants}
			&\frac{\sqrt{|\det D^{2}f(x^{*}(N),N)|}}{\sqrt{|\det D^{2}\widetilde{f}(\widetilde{x}^{*}(N),N)|}}=1+\frac{O(1)}{\sqrt{h(N)}},
		\end{align}
		as $N\to\infty$, where $\widetilde{x}^{*}(N)$ is a maximum of $\widetilde{f}$.
	\end{proposition}
	\begin{proof}
		The proof is analogical to the proof of Proposition 1 in \cite{laplace_paper}.
		The difference is that here the function $\widetilde{f}$ is defined with a more general function $h(N)$ instead of $N$. We can replace $N$ with $h(N)$ everywhere in the results without a significant effect on the proof.
	\end{proof}

	Now, having Proposition \ref{probability_preposition} we can prove the following limit theorems with the estimates valid for a sufficiently large $N$ and the parameter $\delta\in(0,\frac{1}{3(m+1)})$
 
	\begin{theorem}[Central limit theorem I]
		\label{central_limit_theorem_1}
		For $X(N)$ with distribution (\ref{pmf_interior}), the random vector $Z(N)=\sqrt{h(N)}(x^{*}-X(N))$ converges weakly to $\mathcal{N}(0,D^{2}f(x^{*})^{-1})$ and the following estimate of the mgf holds
		\begin{align*}
			M_{Z(N)}(\xi)&=\exp\bigg(\frac{1}{2}\xi^{T}D^{2}f(x^{*})^{-1}\xi\bigg)\times\\
			&\bigg(1+\frac{O(1)}{h(N)^{1/2-3\delta}}+O(1)\frac{h(N)^{1/2+(m+1)\delta}}{N}+O(1)\epsilon(N)\sqrt{h(N)}\bigg).
		\end{align*}
	\end{theorem}

	Here, let us introduce notation $\xi_{y}=(\xi_{2},\ldots,\xi_{m})$, $Y=(X_{2}(N),\ldots,X_{m}(N))$ and $y^{*}=(x^{*}_{2},\ldots,x^{*}_{m})$.

	\begin{theorem}[Central limit theorem II]
		\label{central_limit_theorem_2}
		For $X(N)$ with distribution (\ref{pmf_boundary}) and assuming $\lim_{N\to\infty}\frac{h(N)}{N}=0$, the random vector $Z(N)=\big(h(N)(x_{1}^{*}-X_{1}(N)),\sqrt{h(N)}(y^{*}-Y(N))\big)$ converges weakly to $Exp\big|\frac{\partial f(x^{*})}{\partial x_{1}}\big|$ for $Z_{1}(N)$ and to $\mathcal{N}(0,D_{y}^{2}f(x^{*})^{-1})$ for $\big(Z_{2}(N),\ldots,Z_{m}(N)\big)$. Furthermore, the following estimate of the mgf holds
		\begin{align*}
			M_{Z(N)}(\xi)&=\frac{\big|\frac{\partial f(x^{*})}{\partial x_{1}}\big|}{\big|\frac{\partial f(x^{*})}{\partial x_{1}}-\xi_{1}\big|}\exp\bigg(\frac{1}{2}\xi_{y}^{T}D_{y}^{2}f(x^{*})^{-1}\xi_{y}\bigg)\times\\
			&\bigg(1+\frac{O(1)}{h(N)^{1/2-3\delta}}+O(1)\frac{h(N)^{1/2+(m+1)\delta}}{N}+O(1)\epsilon(N)\sqrt{h(N)}\bigg).
		\end{align*}
	\end{theorem}

	\begin{theorem}[Central limit theorem III]
		\label{central_limit_theorem_3}
		For $X(N)$ with distribution (\ref{pmf_boundary}) and assuming $h(N)=N$ the random vector $Z(N)=\big(h(N)(x_{1}^{*}-X_{1}(N)),\sqrt{h(N)}(y^{*}-Y(N))\big)$ converges weakly to the discrete distribution with the pmf 
		\begin{equation*}
			P(Z_{1}(N)=i)=\exp\bigg(\frac{\partial f(x^{*})}{\partial x_{1}}i\bigg)\bigg(1-\exp\bigg(\frac{\partial f(x^{*})}{\partial x_{1}}\bigg)\bigg),
		\end{equation*}
		for $Z_{1}(N)$ and to $\mathcal{N}(0,D_{y}^{2}f(x^{*})^{-1})$ for $\big(Z_{2}(N),\ldots,Z_{m}(N)\big)$. Furthermore, the following estimate of the mgf holds
		\begin{align*}
			M_{Z(N)}(\xi)&=\frac{1-\exp\big(\frac{\partial f(x^{*})}{\partial x_{1}}\big)}{1-\exp\big(\frac{\partial f(x^{*})}{\partial x_{1}}-\xi_{1}\big)}\exp\bigg(\frac{1}{2}\xi_{y}^{T}D_{y}^{2}f(x^{*})^{-1}\xi_{y}\bigg)\times\\
			&\bigg(1+\frac{O(1)}{h(N)^{1/2-3\delta}}+O(1)\frac{h(N)^{1/2+(m+1)\delta}}{N}+\epsilon(N)\sqrt{h(N)}O(1)\bigg).
		\end{align*}
	\end{theorem}

	\begin{proof}[Proof of Theorem \ref{central_limit_theorem_1}]
		Let us define  $\widetilde{f}(x,N):=f(x,N)+\frac{1}{\sqrt{h(N)}}\xi^{T}(x^{*}-x)$. Then the mgf of $Z(N)$ can be expressed
		\begin{equation*}
			\label{fluctuation_proof_mgf}
			M_{Z(N)}(\xi)=\frac{\sum_{\mathcal{A}\cap L_{N}}e^{h(N)\tilde{f}(x,N)}}{\sum_{\mathcal{A}\cap L_{N}}e^{h(N)f(x,N)}}.
		\end{equation*}
		First, we approximate the numerator and the denominator of the mgf using Theorem \ref{multivariate_interior_theorem}
		\begin{align*}
			M_{Z(N)}(\xi)&=\exp\big(h(N)\widetilde{f}(\widetilde{x}^{*}(N),N)-h(N)f(x^{*}(N),N)\big)\frac{\sqrt{|\det D^{2}f(x^{*}(N),N)|}}{\sqrt{|\det D^{2}\widetilde{f}(\widetilde{x}^{*}(N),N)|}}\times\\
			&\bigg(1+\frac{O(1)}{h(N)^{1/2-3\delta}}+O(1)\frac{h(N)^{1/2+(m+1)\delta}}{N}\bigg).
		\end{align*}
		Next, we insert the estimates (\ref{probability_tilda_estimates_functions}) and (\ref{probability_tilda_estimates_determinants}) from Proposition \ref{probability_preposition}
		\begin{align*}
			M_{Z(N)}(\xi)&=\exp\bigg(\frac{1}{2}\xi^{T}D^{2}f(x^{*})^{-1}\xi+\frac{O(1)}{\sqrt{h(N)}}+O(1)\epsilon(N)\sqrt{h(N)}\bigg)\times\\
			&\bigg(1+\frac{O(1)}{h(N)^{1/2-3\delta}}+O(1)\frac{h(N)^{1/2+(m+1)\delta}}{N}\bigg),
		\end{align*}
		and use the estimate
		\begin{equation*}
			\exp\bigg(\frac{O(1)}{\sqrt{h(N)}}+O(1)\epsilon(N)\sqrt{h(N)}\bigg)=1+\frac{O(1)}{\sqrt{h(N)}}+O(1)\epsilon(N)\sqrt{h(N)},
		\end{equation*}
		which leads to the final result.
	\end{proof}

	\begin{proof}[Proof of Theorem \ref{central_limit_theorem_2}]
		This proof is analogous. Here, we define $\widetilde{f}(x,N):=f(x,N)+\frac{1}{\sqrt{h(N)}}\xi^{T}(\sqrt{h(N)}(x^{*}_{1}-x_{1}),x_{2}^{*}-x_{2},\ldots,x^{*}_{m}-x_{m})$. Then approximate the numerator and the denominator of the mgf using Theorem 3 
		\begin{align*}
			M_{Z(N)}(\xi)&=\exp\big(h(N)\widetilde{f}(\widetilde{x}^{*}(N),N)-Nf(x^{*}(N),N)\big)\frac{1-\big(\frac{h(N)}{N}\frac{\partial f(x^{*}(N),N)}{\partial x_{1}}\big)}{1-\big(\frac{h(N)}{N}\frac{\partial\widetilde{f}(\widetilde{x}^{*}(N),N)}{\partial x_{1}}\big)}\times\\
			&\frac{\sqrt{|\det D_{y}^{2}f(x^{*}(N),N)|}}{\sqrt{|\det D_{y}^{2}\widetilde{f}(\widetilde{x}^{*}(N),N)|}}\bigg(1+\frac{O(1)}{h(N)^{1/2-3\delta}}+O(1)\frac{h(N)^{1/2+(m+1)\delta}}{N}\bigg).
		\end{align*}
		Since the first coordinate of $x^{*}(N)$ is independent of $N$, we can consider $f$ as a function of the remaining $m-1$ coordinates. In that situation, we will have the case (a) of $f$, hence we can apply Proposition \ref{probability_preposition} to obtain
		\begin{align}
			\label{fluctuations_proof_mgf_estimates_middle_step}
			M_{Z(N)}(\xi)&=\frac{1-\exp\big(\frac{h(N)}{N}\frac{\partial f(x^{*}(N),N)}{\partial x_{1}}\big)}{1-\exp\big(\frac{h(N)}{N}\frac{\partial f(\widetilde{x}^{*}(N),N)}{\partial x_{1}}-\frac{h(N)}{N}\xi_{1}\big)}\exp\bigg(\frac{1}{2}\xi_{y}^{T}D^{2}f(x^{*})^{-1}\xi_{y}\bigg)\times\notag\\
			&\bigg(1+\frac{O(1)}{h(N)^{1/2-3\delta}}+O(1)\frac{h(N)^{1/2+(m+1)\delta}}{N}+O(1)\epsilon(N)\sqrt{h(N)}\bigg).
		\end{align}
		Next, using the first order Taylor's Theorem yields
		\begin{equation*}
			\bigg|\frac{\partial f(\widetilde{x}^{*}(N),N)}{\partial x_{1}}-\frac{\partial f(x^{*},N)}{\partial x_{1}}\bigg|\leq F^{(2)}\big(|\widetilde{x}^{*}(N)-x^{*}(N)|+|x^{*}(N)-x^{*}|\big),
		\end{equation*}
		and with use of  (\ref{probability_equation_for_f}), (\ref{probability_estimate_maximums}) and (\ref{probability_tilda_estimates_maximums}) we get
		\begin{equation*}
			\label{probability_fluctuation_proof_b_case_first_derivative_estimate}
			\frac{\partial f(\widetilde{x}^{*}(N),N)}{\partial x_{1}}=\frac{\partial f(x^{*})}{\partial x_{1}}+\frac{O(1)}{\sqrt{h(N)}}.
		\end{equation*}
		Analogous estimation procedure is for the other derivative
		\begin{equation*}
			\label{probability_fluctuation_proof_b_case_first_derivative_estimate}
			\frac{\partial f(x^{*}(N),N)}{\partial x_{1}}=\frac{\partial f(x^{*})}{\partial x_{1}}+O(1)\epsilon(N).
		\end{equation*}
		Using the above estimates and the approximation $e^{x}=1+x+O(x^{2})$ if $x\to 0$, we obtain
		\begin{align*}
			&\Bigg|\frac{1-\exp\big(\frac{h(N)}{N}\frac{\partial f(x^{*}(N),N)}{\partial x_{1}}\big)}{1-\exp\big(\frac{h(N)}{N}\frac{\partial f(\widetilde{x}^{*}(N),N)}{\partial x_{1}}-\frac{h(N)}{N}\xi_{1}\big)}\Bigg|=\Bigg|\frac{\frac{\partial f(x^{*}(N),N)}{\partial x_{1}}+O(1)\frac{h(N)}{N}\big(\frac{\partial f(x^{*}(N),N)}{\partial x_{1}}\big)^{2}}{\frac{\partial f(\widetilde{x}^{*}(N),N)}{\partial x_{1}}-\xi_{1}+O(1)\frac{h(N)}{N}\big(\frac{\partial f(\widetilde{x}^{*}(N),N)}{\partial x_{1}}-\xi_{1}\big)^{2}}\Bigg|=\\
			&\Bigg|\frac{\frac{f(x^{*})}{\partial x_{1}}\big(1+O(1)\epsilon(N)+O(1)\frac{h(N)}{N}\big)}{\big(\frac{\partial f(x^{*})}{\partial x_{1}}-\xi_{1}\big)\big(1+O(1)\frac{h(N)}{N}+\frac{O(1)}{\sqrt{h(N)}}\big)}\Bigg|=\frac{\big|\frac{\partial f(x^{*})}{\partial x_{1}}\big|}{\big|\frac{\partial f(x^{*})}{\partial x_{1}}-\xi_{1}\big|}\bigg(1+\frac{O(1)}{\sqrt{h(N)}}+O(1)\frac{h(N)}{N}\bigg),
		\end{align*}
		and again substituting that into estimate (\ref{fluctuations_proof_mgf_estimates_middle_step}) yields the final result.
	\end{proof}
	
\begin{proof}[Proof of Theorem \ref{central_limit_theorem_3}]
		The proof is similar to the proof of Theorem \ref{central_limit_theorem_2}. Here we have $h(N)=N$ and therefore the main difference is in the last step. That is, using $e^{x}=1+O(x)$ if $x\to 0$ we obtain
		\begin{align*}
			&\Bigg|\frac{1-\exp\big(\frac{\partial f(x^{*}(N),N)}{\partial x_{1}}\big)}{1-\exp\big(\frac{\partial f(\widetilde{x}^{*}(N),N)}{\partial x_{1}}-\xi_{1}\big)}\Bigg|=\frac{\big|1-\exp\big(\frac{\partial f(x^{*})}{\partial x_{1}}\big)\big|}{\big|1-\exp\big(\frac{\partial f(\widetilde{x}^{*})}{\partial x_{1}}-\xi_{1}\big)+\frac{O(1)}{\sqrt{h(N)}}\big|}(1+O(1)\epsilon(N))=\\
			&\frac{\big|1-\exp\big(\frac{\partial f(x^{*})}{\partial x_{1}}\big)\big|}{\big|1-\exp\big(\frac{\partial f(\widetilde{x}^{*})}{\partial x_{1}}-\xi_{1}\big)\big|}\bigg(1+\frac{O(1)}{\sqrt{h(N)}}\bigg),
		\end{align*}
		and substituting that into appropriate estimate of the mgf yields the final result.
	\end{proof}

\begin{acknowledgments}
	Author dedicates special thanks to Professor Vassili Kolokoltsov from University of Warwick, for hints in developing the results of this work. That was, using the Taylor's Theorem in the proofs and, as outlined in the Remark 1, obtaining the solution for some cases by reducing the limit for natural numbers to the limit for specific subsequence of the natural numbers.
\end{acknowledgments}


\bibliographystyle{apa}
\bibliography{biblio}
	
\end{document}